\documentclass[a4paper]{amsart}

\usepackage{amsmath,amsthm,amsfonts,amssymb,amscd}
\input{xy}
\xyoption{all}

\newtheorem{theo}{Theorem}[section]
\newtheorem{pro}[theo]{Proposition}
\newtheorem{lemma}[theo]{Lemma}
\newtheorem{coro}[theo]{Corollary}
\newtheorem{defi}[theo]{Definition}

\newcommand{\dem}{{\textit{Proof. }}}

\newenvironment{sis}{\left\{\begin{aligned}}{\end{aligned}\right.}

\numberwithin{equation}{section}

\newcommand{\Z}{\mathbb{Z}}
\newcommand{\N}{\mathbb{N}}

\newcommand{\An}{A(n)=F[x_1,\dots , x_n]/(x_1^p, \dots, x_n^p)}

\newcommand{\0}{\underline{0}}

\newcommand{\g}{\mathfrak{g}}
\newcommand{\h}{\mathfrak{h}}

\newcommand{\ep}{\epsilon}

\renewcommand{\d}{{\rm d}}

\newcommand{\Sq}{{\rm Sq}}

\begin{document}


\title{Infinitesimal deformations of restricted simple Lie algebras I}
\author{Filippo Viviani}

\date{2 July 2008}

\address{Institut f\"ur Mathematik, Humboldt Universit\"at zu Berlin, 10099 Berlin (Germany).}
\email{viviani@amath.hu-berlin.de}



\thanks{During the preparation of this paper, the author was partially supported 
by a grant from the Mittag-Leffler Institute in Stockholm.}

\maketitle

\begin{abstract}
We compute the infinitesimal deformations of two families of restricted simple
modular Lie algebras of Cartan-type: the Witt-Jacobson
and the Special Lie algebras.
\end{abstract}



\section{Introduction}

Simple Lie algebras over an algebraically closed field of
characteristic zero were classified at the beginning of the XIX
century by Killing and Cartan. They used the non-degeneracy of the
Killing form to describe the simple Lie algebras in terms of root
systems which are then classified by Dynkin diagrams.

 This method breaks down in positive characteristic because the Killing form may
degenerate. Indeed the classification problem remained open for a
long time until it was recently solved, if the characteristic of the
base field is greater than $3$, by Wilson-Block (see \cite{BW}),
Strade (see \cite{STR1}, \cite{STR2}, \cite{STR3}, \cite{STR4},
\cite{STR5}, \cite{STR6}) and Premet-Strade (see \cite{PS1},
\cite{PS2}, \cite{PS3}). The classification remains still open in
characteristic $2$ and $3$ (see \cite[page 209]{STR}).

According to this classification, \emph{simple modular} (that is over a field of positive
characteristic) Lie algebras are divided into two big families, called classical-type
and Cartan-type algebras.
The algebras of classical-type are obtained by the simple Lie algebras
in characteristic zero by first taking a model over the integers (via Chevalley bases)
and then reducing modulo $p$ (see \cite{SEL}). The algebras of Cartan-type were constructed by
Kostrikin-Shafarevich in 1966 (see \cite{KS}) as finite-dimensional analogues
of the infinite-dimensional complex simple Lie algebras, which occurred in Cartan's classification
of Lie pseudogroups, and are divided into four families, called Witt-Jacobson, Special,
Hamiltonian and Contact algebras. The Witt-Jacobson Lie algebras
are derivation algebras of truncated divided power algebras and the remaining three families
are the subalgebras of derivations fixing a volume form, a Hamiltonian form and a contact form,
respectively. Moreover in characteristic $5$ there is one exceptional simple modular Lie algebra
called the Melikian algebra (introduced in \cite{MEL}).

We are interested in a particular class of modular Lie algebras
called \emph{restricted}. These can be characterized as those
modular Lie algebras such that the $p$-power of an inner derivation
(which in characteristic $p$ is a derivation) is still inner.
Important examples of restricted Lie algebras are the ones coming
from groups schemes. Indeed there is a one-to-one correspondence
between restricted Lie algebras and finite group schemes whose
Frobenius vanishes (see \cite[Chap. 2]{DG}).

By standard facts of deformation theory, the \emph{infinitesimal deformations} of a
Lie algebra are parametrized by
the second cohomology of the Lie algebra with values in the adjoint representation
(see for example \cite{GER1}).

It is a classical result (see \cite{HiSt}) that for a simple Lie
algebra $\g$ over a field of characteristic $0$ it holds that
$H^i(\g,\g)=0$ for every $i\geq 0$, which implies in particular that
such Lie algebras are rigid. The proof of this fact relies on the
non-degeneracy of the Killing form and the non-vanishing of the
trace of the Casimir element, which is equal to the dimension of the
Lie algebra. Therefore the same proof works also for the simple
modular Lie algebras of classical type over a field of
characteristic not dividing the determinant of the Killing form and
the dimension of the Lie algebra. Actually Rudakov (see \cite{RUD})
showed that such Lie algebras are rigid if the characteristic of the
base field is greater than or equal to $5$ while in characteristic $2$
and $3$ there are non-rigid classical Lie algebras (see \cite{Che1},
\cite{Che2}, \cite{Che3}).

The purpose of this article is to compute the infinitesimal deformations of the
first two families of restricted simple
Lie algebras of Cartan type: the Witt-Jacobson algebras $W(n)$ and the Special
algebras $S(n)$. Unlike the classical-type simple algebras, it turns out that
these two families are not rigid. More precisely we get the
following two Theorems (we refer to
subsections $3.1$ and $4.1$ for the standard notations concerning $W(n)$ and $S(n)$
and to subsection $2.3$ for the definition of the squaring operators $\Sq$).

\begin{theo}\label{W-finaltheorem}
Assume that the characteristic $p$ of the base field $F$ is different from $2$. Then
we have
$$H^2(W(n),W(n))=\bigoplus_{i=1}^n F\cdot \langle {\rm Sq}(D_i)\rangle$$
with the exception of the case $n=1$ and $p=3$ when it is $0$.
\end{theo}

\begin{theo}\label{S-finaltheorem}
Assume that the characteristic of the base field $F$ is different from $2$ and moreover
it is different from $3$ if $n=3$. Then
we have
$$H^2(S(n),S(n))=\bigoplus_{i=1}^n F\cdot \langle {\rm Sq}(D_i)\rangle   \bigoplus
F\cdot \langle\Theta \rangle$$
where $\Theta$ is defined by $\Theta(D_i,D_j)=D_{ij}(x^{\tau})$ and extended by $0$ outside  $S(n)_{-1}\times S(n)_{-1}$.
\end{theo}

In the two forthcoming papers \cite{VIV2, VIV3}, we compute the infinitesimal deformations of
the remaining restricted simple Lie algebras of Cartan-type, namely the
Hamiltonian, the Contact and the exceptional Melikian algebras. Moreover, in another
paper \cite{VIV4}, we apply these results to the study of the infinitesimal 
deformations of the simple finite group schemes corresponding to the restricted simple 
Lie algebras of Cartan type.

Let us mention that the infinitesimal deformations of simple Lie
algebras of Cartan-type (in the general non-restricted case) have
been considered already by D\v zumadildaev in \cite{DZ1, DZ2, DZ3} and 
D\v zumadildaev-Kostrikin in \cite{DK} but a
complete picture as well as detailed proofs were missing. More
precisely: in \cite{DK} the authors compute the infinitesimal
deformations of the Jacobson-Witt algebras of rank $1$, in
\cite[Theorem 4]{DZ1} the author describes the infinitesimal
deformation of the Jacobson-Witt algebras of any rank but without a
proof, in \cite{DZ2} a general strategy for the Jacobson-Witt and
Hamiltonian algebras is outlined (without proofs) and finally in
\cite{DZ3} the author clarifies this strategy and then applies it to
the Jacobson-Witt algebras but with a half-page sketch of the proof.

Our approach works for all the restricted simple Lie algebras of Cartan-type  and
is different from the approach of D\v zumadildaev
although we took from him the idea to consider relative cohomology with respect to
the subalgebra of negative degree elements.
As a byproduct of our proof, we recover the results of Celousov (see \cite{CEL})
on the first cohomology group of the adjoint representation (Theorems \ref{W-Celousov}
and \ref{S-Celousov}).

The results presented here constitute part of my doctoral thesis. I
thank my advisor prof. R. Schoof for useful advice and constant
encouragement.

\section{Some preliminaries results on the cohomology of Lie algebras}

\subsection{Review of general theory}

In this subsection we review, in order to fix notations, the
classical theory of cohomology of Lie algebras (see for example
\cite{HS}).

If $\g$ is a Lie algebra over a field $F$ and $M$ is a $\g$-module, then
the cohomology groups $H^*(\g,M)$ can be computed from the
complex of $n$-dimensional cochains $C^n(\g, M)$ ($n\geq 0$),
that are alternating $n$-linear functions $f:\Lambda^{n}(\g)\to M$,
with differential $d:C^n(\g,M)\to C^{n+1}(\g,M)$ defined by
\begin{equation} \label{d} \begin{split}
\d f(\sigma_0,\dots,\sigma_n)=&\sum_{i=0}^n(-1)^i\sigma_i\cdot f(\sigma_0,\dots, \hat{\sigma_i},
\dots, \sigma_n)+\\
&\sum_{p<q}(-1)^{p+q}f([\sigma_p,\sigma_q],\sigma_0,\dots,\hat{\sigma_p},\dots, \hat{\sigma_q},
\dots \sigma_n),
\end{split} \end{equation}
where the sign $\hat{}$ means that the argument below must be omitted. Given 
$f \in C^n(\g, M)$ and $\gamma \in \g$, we denote with $f_{\gamma}$ the restriction
of $f$ to $\gamma\in \g$, that is the element of $C^{n-1}(\g, M)$ given by 
$$f_{\gamma}(\sigma_0, \cdots, \sigma_{n-1}):=f(\gamma, \sigma_0, \cdots, 
\sigma_{n-1}).$$
With this notation, the above differential satisfies the 
following useful formula (for any $\gamma\in \g$ and $f\in C^n(\g,M)$):
\begin{gather}
\d(\gamma \cdot f)=\gamma\cdot (\d f)\label{module-d},\\
(\d f)_{\gamma}=\gamma\cdot f-\d(f_{\gamma}).\label{restriction}
\end{gather}
where each $C^n(\g,M)$ is a $\g$-module by means of the action
\begin{equation}\label{module}
(\gamma\cdot f)(\sigma_1,\dots,\sigma_n)=\gamma\cdot f(\sigma_1,\dots, \sigma_n)
-\sum_{i=1}^n f(\sigma_1,\cdots,[\gamma,\sigma_i],\dots \sigma_n).
\end{equation}

As usual we indicate with $Z^n(\g,M)$ the subspace of $n$-cocycles
and with $B^n(\g,$ $ M)$ the subspace of $n$-coboundaries. Therefore
$H^n(\g,M):=Z^n(\g,M)/$ $B^n(\g,M)$.

A useful tool to compute cohomology of Lie algebras is the following
Hochschild-Serre spectral sequence relative to a subalgebra $\h<\g$:

\begin{equation}\label{HS-ss1}
E_1^{p,q}=H^q(\h, C^p(\g/\h, M))\Longrightarrow H^{p+q}(\g, M),
\end{equation}
which in the case where $\h$ is an ideal of $\g$ (which we indicate as $\h\lhd \g$) becomes
\begin{equation}\label{HS-ss2}
E_2^{p,q}=H^p(\g/\h, H^q(\h, M))\Longrightarrow H^{p+q}(\g, M).
\end{equation}
Moreover for the second page of the first spectral sequence
\ref{HS-ss1}, we have the equality
\begin{equation}\label{REL-coho}
E_2^{p,0}= H^p(\g,\h;M),
\end{equation}
where $H^*(\g,\h;M)$ are the relative cohomology groups defined (by
Chevalley and Eilenberg \cite{CE}) from the sub-complex
$C^p(\g,\h;M)\subset C^p(\g,M)$ consisting of cochains orthogonal to
$\h$, that is cochains satisfying the two conditions:
\begin{gather}
f_{| \h}=0, \label{ortho1}\\
\d f_{| \h}=0 \hspace{0,5cm} \text{ or equivalently } \hspace{0,5cm}
\gamma \cdot f=0 \hspace{0,3cm} \text{ for every }
\gamma\in \h.\label{ortho2}
\end{gather}
Note that in the case where $\h\lhd \g$, the equality \ref{REL-coho}
is consistent with the second spectral sequence \ref{HS-ss2} because in that case
we have $H^p(\g, \h, M)=H^p(\g/\h, M^{\h})$.

\subsection{Torus actions and Gradings}

The Lie algebras that we consider in this paper, namely the Witt-Jacobson
Lie algebra $W(n)$ and the Special algebra $S(n)$, are graded algebras
which admit a root space decomposition with respect to a maximal torus
contained in the $0$-graded piece. Under these hypothesis, the cohomology groups
admit a very useful decomposition that we are going to review in this subsection.

Suppose that a torus $T$ acts on both $\g$ and $M$ in a way that is
compatible with the action of $\g$ on $M$, which means that $t\cdot(g\cdot m)=(t\cdot g)\cdot m+
t\cdot (g\cdot m)$ for every $t\in T$, $g\in \g$ and $m\in M$. Then the
action of $T$ can be extended to the space of $n$-cochains by
$$(t\cdot f)(\sigma_1,\cdots ,\sigma_n)=t\cdot f(\sigma_1, \cdots, \sigma_n)
-\sum_{i=1}^n f(\sigma_1, \cdots, t\cdot \sigma_i, \cdots, \sigma_n).$$
It follows easily from the compatibility of the action of $T$ and formula
\ref{restriction}, that the action of $T$ on the cochains commutes with the
differential $\d$. Therefore, since the action of a torus is always completely
reducible, we get a decomposition in eigenspaces
\begin{equation}\label{deco1}
H^n(\g,M)=\bigoplus_{\phi\in \Phi} H^n(\g,M)_{\phi},
\end{equation}
where $\Phi={\rm Hom}_F(T,F)$ and $H^n(\g,M)_{\phi}=\{[f] \in
H^n(\g,M)\: | \: t\cdot [f]=\phi(t) [f] \: \: \text{ if }$ $  t\in
T\}$. A particular case of this situation occurs when $T\subset \g$
and $T$ acts on $\g$ via the adjoint action and on $M$ via
restriction of the action of $\g$. It is clear that this action is
compatible and moreover the above decomposition reduces to
$$H^n(\g,M)=H^n(\g,M)_{\0} \: , $$
where $\underline{0}$ is the trivial homomorphism (in this situation we say that
the cohomology reduces to \emph{homogeneous} cohomology). Indeed, if we consider
an element $f\in Z^n(\g,M)_{\phi}$, then by applying formula \ref{restriction}
with $\gamma=t\in T$ we get
$$0=(\d f)_t=t\cdot f- \d (f_t)=\phi(t)f-\d(f_t),$$
from which we see that the existence of a $t\in T$ such that $\phi(t)\neq 0$
forces $f$ to be a coboundary.

Now suppose that $\g$ and $M$ are graded and that the action of $\g$ respects
these gradings, which means that $\g_d\cdot M_e\subset M_{d+e}$ for all $e, d\geq 0$.
Then the space
of cochains can also be graded: a homogeneous cochain $f$ of degree $d$ is
a cochain such that $f(\g_{e_1}\times \cdots \times \g_{e_n})\subset M_{\sum e_i+d}$.
With this definition, the differential becomes of degree $0$ and therefore we get
a degree decomposition
\begin{equation}\label{deco2}
H^n(\g,M)=\bigoplus_{d\in \Z} H^n(\g,M)_d.
\end{equation}
Finally, if the action of $T$ is compatible with the grading, in the sense
that $T$ acts via degree $0$ operators both on $\g$ and on $M$, then the above
two decompositions \ref{deco1} and \ref{deco2} are compatible and give rise to the refined
weight-degree decomposition
\begin{equation}\label{deco3}
H^n(\g,M)=\bigoplus_{\phi\in \Phi}\bigoplus_{d\in \Z}H^n(\g,M)_{\phi,d}.
\end{equation}

\subsection{Squaring operation}

There is a canonical way to produce $2$-cocycles in $Z^2(\g,\g)$ over a field of characteristic
$p>0$, namely the squaring operation (see \cite{GER1}).
Given a derivation $\gamma\in Z^1(\g,\g)$ (inner or not), one defines the squaring of
$\gamma$ to be
\begin{equation}\label{Square}
\Sq(\gamma)(x,y)=\sum_{i=1}^{p-1}\frac{[\gamma^i(x),\gamma^{p-i}(y)]}{i!(p-i)!}\in Z^2(\g,\g),
\end{equation}
where $\gamma^i$ is the $i$-th iteration of $\gamma$.
In \cite{GER1} it is shown that $[\Sq(\gamma)]\in H^2(\g,\g)$ is an obstruction to integrability
of the derivation $\gamma$, that is to the possibility of finding an automorphism
of $\g$ extending the infinitesimal automorphism given by $\gamma$.

\section{The Witt-Jacobson algebra}

\subsection{Definition and basic properties}

We first introduce some useful notations.
Inside the set $\Z^n$ of $n$-tuples of integers,
we consider the order relation defined by
$a=(a_1, \dots, a_n)< b=(b_1, \dots, b_n)$ if $a_i<b_i$ for every $i=1, \dots, n$.
We call degree of $a\in \Z^n$ the number $|a|=\sum_{i=1}^n a_i$. For every integer
$0\leq l<p$, we define $\underline{l}:=(l,\cdots, l)$ and we set
$\tau:=\underline{p-1}$ (this $n$-tuple will appear often in what follows and hence it deserves
a special notation). Moreover, for every $j\in \{1,\dots, n\}$ we call $\epsilon_j$ the
$n$-tuple having $1$ at the $j$-th place and $0$ otherwise.

Let $A(n)=F[x_1,\dots , x_n]/(x_1^p, \dots, x_n^p)$ be the ring of $p$-truncated polynomial
in $n$ variables over a field $F$ of positive characteristic $p>0$.
Note that $A(n)$ is a finite $F$-algebra of dimension $p^n$ with a basis given by the elements
$\{x^a:=x_1^{a_1}\dots x_n^{a_n} \: |\: a\in \Z^n,  \: \underline{0}\leq a \leq \tau\}$.
Moreover it has a natural graduation $A(n)=\oplus_{i=0}^{n(p-1)} A(n)_i$,
obtained by assigning to the monomial $x^a$ the degree $|a|$.

\begin{defi}
The Witt-Jacobson algebra $W(n)$ is the restricted Lie algebra ${\rm Der}_FA(n)$
of derivations of $\An$.
\end{defi}

For every $j\in \{1, \dots , n\}$, we put $D_j:=\frac{\partial}{\partial x_j}$.
The Witt-Jacobson algebra $W(n)$ is
a free $A(n)$-module with basis $\{D_1, \dots, D_n\}$. Hence ${\rm dim}_F(W(n))=np^n$ with a
basis over $F$ given by $\{x^aD_j\: | \: 1\leq j\leq n, \: \underline{0}\leq a \leq \tau\}$.

Moreover $W(n)$ is a graded Lie algebra with the $\Z$-gradation defined by
$W(n)_i:=\sum_{j=1}^n A(n)_{i+1}D_j$
where $i=-1, \dots, n(p-1)-1$.
Note that the unique summand of negative degree is $W(n)_{-1}=\oplus_{i=1}^n F\cdot
\langle D_i \rangle$ while
the summand of degree $0$ is $W(n)_0=\oplus_{1\leq i, j\leq n}F\cdot \langle x_iD_j \rangle$
and its adjoint action on $W(n)_{-1}$ induces an isomorphism $W(n)_0\cong  \mathfrak{gl}(n, F)$.

The algebra $W(n)$ is simple unless $p=2$ and $n=1$ (see \cite[Chap. 4, Theo. 2.4]{FS})
and it admits a root space decomposition with respect to a canonical Cartan subalgebra.

\begin{pro} \label{Cartan-W}
For each $i\in \{1, \dots, n\}$,  let $h_i=x_iD_i$.
\begin{itemize}
\item[(a)] $T:=\sum_{i=1}^n Fh_i$ is a maximal torus of $W(n)$ (called the canonical maximal
torus).
\item[(b)] The centralizer of $T$ inside $W(n)$ is $T$ itself, which is hence a Cartan
subalgebra of $W(n)$.
\item[(c)] Let $\Phi:={\rm Hom}_{\mathbb{F}_p}(\oplus_{i=1}^n\mathbb{F}_p \cdot h_i,
\mathbb{F}_p)$, where $\mathbb{F}_p$ is the prime field
of $F$. In the Cartan decomposition
$W(n)=\oplus_{\phi\in \Phi} W(n)_{\phi}$, every direct summand $W(n)_{\phi}$
has dimension $n$. Moreover $x^aD_i\in W(n)_{a-\epsilon_i}$, where $a-\epsilon_i$
is viewed as an element of $\Phi$ by reduction modulo $p$.
\end{itemize}
\end{pro}
\dem
 See \cite[Chap. 4, Theo. 2.5]{FS}.
\hfill \qed

\subsection{Strategy of the proof of the Main Theorem}

In this subsection we outline the strategy of the proof of Theorem
\ref{W-finaltheorem} from the Introduction. In particular, from now on, we assume that
the base field $F$ has characteristic $p\geq 3$.
Note that in the exceptional case $n=1$ and $p=3$, one has the
isomorphism $W(1)\cong \mathfrak{sl}_2$ and hence we recover the known vanishing
result for the simple algebras of classical-type.

We first observe that the $2$-cocycles $\Sq(D_i)$ appearing in Theorem
\ref{W-finaltheorem} are independent modulo coboundaries unless $n=1$ and $p=3$, in which
case it is easily seen that $\Sq(D_1)=0$.
Indeed, on one hand, for every $g\in C^1\left(W(n),W(n)\right)$ and
$1\leq r,s\leq n$, the following element
$$\d g(x_r^2D_s,x_r^{p-2}x_sD_s)=[x_r^2D_s,g(x_r^{p-2}x_sD_s)]-[x_r^{p-2}x_sD_s,g(x_r^2D_s)]$$
cannot contain terms of negative degree. On the other hand, we get that
\begin{equation}\label{W-Sq}
{\rm Sq}(D_i)(x_r^2D_s,x_r^{p-2}x_sD_s)=
\begin{sis}
&D_s & \text{ if } i=r\neq s,\\
&-3D_i & \text{ if } i=r=s,\\
& 0 & \text{ otherwise} ,
\end{sis}
\end{equation}
which shows the independence of the ${\rm Sq}(D_i)$ modulo coboundaries, using the first
case if $n\geq 2$ and the second if $p\geq 5$.

The proof that these $2$-cocycles generate the whole second cohomology group
is divided into three steps.

\underline{STEP I}: We prove that we can reduce to relative cohomology (see section $2.1$)
with respect to the subalgebra $W(n)_{-1}$ of negative terms:
\begin{equation*}
H^2(W(n),W(n))=H^2(W(n),W(n)_{-1};W(n)).
\end{equation*}
This is achieved by first observing that the second cohomology groups reduces
to homogeneous cohomology with respect to the maximal torus $T<W(n)$ (see section $2.2$)
and then by considering the homogeneous Hochschild-Serre
spectral sequence associated to the subalgebra $W(n)_{-1}< W(n)$ (see (\ref{HS-ss1})):
\begin{equation}\label{W-HS1}
(E_1^{r,s})_{\0}=H^s(W(n)_{-1},C^r(W(n)/W(n)_{-1},W(n)))_{\0}\Rightarrow
H^{r+s}(W(n),W(n))_{\0}.
\end{equation}
We prove that $(E_1^{0,1})_{\0}=(E_1^{0,2})_{\0}=0$
(corollary \ref{W(n)_{-1}-vanishing})
and $(E_2^{1,1})_{\0}=0$ (proposition \ref{WE_2^{1,1}}) which gives the conclusion
by (\ref{REL-coho}).

\underline{STEP II}: Using orthogonality with respect to $W(n)_{-1}$ (see
(\ref{ortho1}) and (\ref{ortho2})), we prove in proposition \ref{W-rel-coho} that
$$H^2(W(n),W(n)_{-1};W(n))= H^2(W(n)_{\geq 0},W(n)_{-1})$$
where $W(n)_{\geq 0}$ acts on $W(n)_{-1}$ by the projection onto $W(n)_{\geq 0}
/W(n)_{\geq 1}=W(n)_0$ followed by the adjoint representation of $W(n)_0=
\mathfrak{gl}(n,F)$ on $W(n)_{-1}$.

Then, by using the Hochschild-Serre spectral sequence with respect to the
ideal $W(n)_{\geq 1}\lhd W(n)_{\geq 0}$ (see (\ref{HS-ss2})),
we prove in proposition \ref{W-H^1} that
$$H^2(W(n)_{\geq 0},W(n)_{-1})=H^2(W(n)_{\geq 1},W(n)_{-1})^{W(n)_0}$$
where $W(n)_{-1}$ is considered as a trivial $W(n)_{\geq 1}$-module.

\underline{STEP III}: We compute the invariant second cohomology group 
$H^2(W(n)_{\geq 1},$ $ W(n)_{-1})^{W(n)_0}$ showing that (unless $p=3$
and $n=1$) it is generated by the
projection onto $W(n)_{-1}$ of the cocycles $\Sq(D_i)$ (proposition
\ref{W-H^2}). The idea of the proof is to approximate this cohomology
group by the truncated cohomology groups
$$H^2\left(\frac{W(n)_{\geq 1}}{W(n)_{\geq d}},W(n)_{-1}\right)^{W(n)_0}$$
which for large $d$ are equal to our cohomology group. The
computation proceeds by induction on $d$ using the Hochschild-Serre
spectral sequence with respect to the ideal
$$\frac{W(n)_{\geq d}}{W(n)_{\geq d+1}}\lhd
\frac{W(n)_{\geq 1}}{W(n)_{\geq d+1}}.$$
In the course of the proof of the main Theorem,
we obtain a new proof of the following result.

\begin{theo}[Celousov]\label{W-Celousov}
$H^1(W(n),W(n))=0.$
\end{theo}
\dem The proof follows the same steps as in the proof of the main
Theorem. The spectral sequence (\ref{W-HS1}), in view of the
corollary \ref{W(n)_{-1}-vanishing} and the Formula
(\ref{REL-coho}), gives that
$$H^1(W(n),W(n))=H^1(W(n),W(n)_{-1};W(n)).$$
Then the required vanishing follows from propositions
\ref{W-rel-coho} and \ref{W-H^1}. \qed

\subsection{Reduction to $W(n)_{-1}$-relative cohomology}

This subsection is devoted to the first step of the proof (see subsection $3.2$),
namely the reduction to the relative cohomology with respect to the subalgebra
$W(n)_{-1}<W(n)$.
First of all we want to prove the vanishing of the homogeneous cohomology groups
$H^s(W(n)_{-1}, W(n))_{\0}$
appearing in the first column of the spectral sequence (\ref{W-HS1}).
For that purpose, we need the following proposition, in which the action of
$W(n)_{-1}$ on $A(n)$ is the natural one.

\begin{pro} \label{W(n)_{-1}-cohomology}
For every $i=1,\cdots,n$, we denote with $x_i^{p-1}D_i^*$ the linear function
from $W(n)_{-1}$ to $A(n)$ which sends $D_i$ to $x_i^{p-1}$ and $D_j$ to $0$ for
$j\neq i$. Then we have
$H^s(W(n)_{-1}, A(n))=\bigwedge^s \oplus_{i=1}^n F \cdot \langle x_i^{p-1}
D_i^* \rangle.$
\end{pro}
\dem Clearly the cochains appearing in the statement are cocycles
and they are independent modulo coboundaries since it follows
easily, from formula (\ref{d}), that if $g\in C^{s-1}(W(n)_{-1},
A(n))$ then $\d g(D_{i_1},\cdots, D_{i_s})\in A(n)$ cannot contain
the monomial $x_{i_1}^{p-1}\cdots x_{i_s}^{p-1}$.

In order to prove that the above cocycles generate the whole cohomology group,
we proceed by double induction
on $s$ and $n$, the case $s=0$ being true since $A(n)^{W(n)_{-1}}=F\cdot 1$.
We view $A(n-1)$  inside $A(n)$ as the subalgebra of polynomials in the variables
$x_2, \cdots, x_n$ and $W(n-1)_{-1}$ inside $W(n)_{-1}$ as the subalgebra generated
by $D_2,\cdots, D_n$. Thus the action of $W(n)$ on $A(n)$ restricts to the
natural action of $W(n-1)$ on $A(n-1)$.

Consider $f\in Z^s(W(n)_{-1}, A(n))$. By adding a coboundary
$\d g$ and using formula (\ref{restriction}) for $\d g$ and $\gamma=D_1$,
we can suppose that
$$f_{| D_1}: W(n-1)_{-1}^{s-1}\to x_1^{p-1} A(n-1).$$
Moreover, since $f$ is a cocycle, the same formula (\ref{restriction}) gives
$$0=(\d f)_{D_1}=[D_1,f(-)]-\d(f_{| D_1}).$$
Now observe that, by the condition above, $\d(f_{| D_1})$ takes values in $x_1^{p-1}A(n-1)$
while obviously $[D_1, f(-)]$ cannot contain monomials with the $x_1$ erased to the
$(p-1)$-th power. Hence it follows that
$$\begin{sis}
&[D_1, f_{| W(n-1)_{-1}^s}]=0,\\
&\d(f_{| D_1})=0.\\
\end{sis}$$
The first equation says that $f_{| W(n-1)_{-1}^s}$ takes values in
$A(n-1)$ and hence belongs to $Z^s(W(n-1)_{-1}, A(n-1))$. The second
equation says that $f_{| D_1}\in
Z^{s-1}(W(n-1)_{-1},x_1^{p-1}A(n-1))=Z^{s-1}(W(n-1)_{-1},A(n-1))
\otimes <x_1^{p-1}>$. In both cases, by induction, we get that $f\in
B^s(W(n)_{-1}, A(n)) + \bigwedge^s \oplus_{i=1}^n F\cdot \langle
x_i^{p-1}D_i^* \rangle$ and this concludes the proof. \qed

\begin{coro}\label{W(n)_{-1}-vanishing}
We have $H^s(W(n)_{-1},W(n))\cong H^s(W(n)_{-1},A(n))\otimes W(n)_{-1}$. \\
Therefore  $(E_1^{0,s})_{\0}=H^s(W(n)_{-1},W(n))_{\0}=0$ for every $s\geq 0$.
\end{coro}
\dem The first claim follows from the $W(n)_{-1}$-decomposition
$W(n)=A(n)\otimes W(n)_{-1}$ and the fact that $W(n)_{-1}$ is an
abelian Lie algebra. The second claim follows from the first and the
fact that $H^s(W(n)_{-1},A(n))=H^s(W(n)_{-1},A(n))_{\0}$ (by
proposition \ref{W(n)_{-1}-cohomology}) while $(W(n)_{-1})_{\0}=0$.
\qed

Now we deal with the term in position $(1,1)$ of the above spectral sequence.
We prove that it vanishes starting from the second level.

\begin{pro}\label{WE_2^{1,1}}
In the spectral sequence (\ref{W-HS1}), we have that
$(E_2^{1,1})_{\0}=0.$
\end{pro}
\dem We have to show the injectivity of the level $1$ differential
map
$$d:(E_1^{1,1})_{\0}\longrightarrow(E_1^{2,1})_{\0}.$$
In the course of this proof, we adopt the following convention: given an element
$f\in  C^1(W(n)_{-1}, C^s(W(n)/W(n)_{-1},W(n))$,
 we write its value on $D_i\in W(n)_{-1}$ as $f_{D_i}\in C^s(W(n)/W(n)_{-1},W(n))$.

We want to show, by induction on the degree of $E\in W(n)/W(n)_{-1}$,
that if $[df]=0\in H^1(W(n)_{-1},C^2(W(n)/W(n)_{-1},W(n)))$ then we can choose a representative
$\widetilde{f}$ of $[f]\in H^1(W(n)_{-1},C^1(W(n)/W(n)_{-1},W(n)))$ such that
$\widetilde{f}_{D_i}(E)=0$ for every $i=1, \cdots, n$. So suppose that we have already found
a representative $f$
 such that $f_{D_i}(F)=0$ for every $F\in W(n)/W(n)_{-1}$ of degree less than $d$
and for every $i$.
First of all, we can find a representative $\widetilde{f}$ of $[f]$ such that
\begin{equation*}\tag{*}
\widetilde{f}_{D_i}(E)\in \langle x_i^{p-1} \rangle \otimes W(n)_{-1}
\end{equation*}
for every $i$ and for every $E\in W(n)$ of degree $d$. Indeed, by
the induction hypothesis, the cocycle condition for $f$ is $\partial
f_{D_i,D_j}(E)=[D_i,f_{D_j}(E)]- [D_j,f_{D_i}(E)]$. On the other
hand, by choosing an element $h\in C^1(W(n)/$ $W(n)_{-1},W(n))$ that
vanishes on the elements of degree less than $d$, we can add to $f$
(without changing its  cohomological class neither affecting the
inductive assumption) the coboundary $\partial h$ whose value on $E$
is $\partial h_{D_i}(E)=[D_i,h(E)]$.  Hence, for a fixed element
$E$ of degree $d$, the map $D_i\mapsto f_{D_i}(E)$ gives rise to an
element of $H^1(W(n)_{-1},W(n))$ and, by proposition
\ref{W(n)_{-1}-cohomology}, we can chose an element $h(E)$ as above
such that the new cochain $\widetilde{f}=f+\partial h$ verifies the
condition  $(*)$ as above.

Note that, by the homogeneity of our cocycles, the functions $\widetilde{f}_{D_i}$ can assume
non-zero values only on the elements $E$ of weights $-\epsilon_k$, for a certain $k$, which are
the form $E=x_k^{p-1}x_hD_h$ for some $k\neq h$ (note that we have already done in the case
$n=1$). Hence, from now on, we can assume that $d=p-1\geq 2$ and pay attention only
to the elements of the above form.

Now we are going to use the condition that $[d\widetilde{f}]=0\in (E_1^{2,1})_{\0}$,
that is $d\widetilde{f}=\partial g$ for some $g\in C^2(W(n)/W(n)_{-1},W(n))_{\0}$.
Explicitly, for $A, B\in W(n)/W(n)_{-1}$) we have that
\begin{equation}\label{Wfor-g}
\partial g_{D_i}(A,B)=[D_i,g(A,B)]-g([D_i,A],B)-g(A,[D_i,B]),
\end{equation}
\begin{equation}\label{Wfor-d}
\begin{split}
d\widetilde{f}_{D_i}(A,B)=
&\widetilde{f}_{D_i}([A,B])-[A,\widetilde{f}_{D_i}(B)]+[B,\widetilde{f}_{D_i}(A)]+\\
&-\delta_{{\rm deg}(A),0}\widetilde{f}_{[D_i,A]}(B)+\delta_{{\rm deg}(B),0}
\widetilde{f}_{[D_i,B]}(A),
\end{split}
\end{equation}
where the last two terms in the second formula are non-zero only if ${\rm deg}(A)=0$ and
${\rm deg}(B)=0$ respectively.
We apply the above formulas for the elements $A=x_k^{p-2}x_h^2D_h$ and $B=x_kD_h$.
Taking into account the inductive hypothesis on the degree and the homogeneity assumptions,
formula (\ref{Wfor-d}) becomes
$$d\widetilde{f}_{D_i}(x_k^{p-2}x_h^2D_h,x_k D_h)=
-2\widetilde{f}_{D_i}(x_k^{p-1}x_hD_h)
=\alpha x_i^{p-1}D_k$$
for a certain $\alpha\in F$, while formula (\ref{Wfor-g}) gives
$$\partial g_{D_i}(x_k^{p-2}x_h^2D_h,x_kD_h)=[D_i,g(x_k^{p-2}x_h^2D_h,x_kD_h)]-
g([D_i,x_k^{p-2}x_h^2D_h],x_kD_h).$$
Observe that if ${\rm deg}(B)=0$ and ${\rm deg}(A)<p-1$, then
${\rm deg}_{x_i}(g(A,B))\leq $ ${\rm deg}_{x_i}(A)$ (where ${\rm deg}_{x_i}(-)$ indicate the largest
power of $x_i$ which appears in the argument). Indeed, by the inductive hypothesis, formula
(\ref{Wfor-d}) gives that $d\widetilde{f}_{D_i}(A,B)=0$ and hence the conclusion follows by
repeatedly applying formula (\ref{Wfor-g}): $0=\partial g_{D_i}(A,B)=[D_i,g(A,B)]-g([D_i,A],B)$.

From this observation, it follows that
$g([D_i,x_k^{p-2}x_h^2D_h],x_kD_h)$ cannot contain a monomial of the
form $x_i^{p-1}D_k$ and hence neither can the element \\ $\partial
g_{D_i}(x_k^{p-2}x_h^2D_h,x_kD_h)$, since in the above formula the
first element is a derivation with respect to $D_i$. Therefore by
imposing $d\widetilde{f}_{D_i}=\partial g_{D_i}$, we obtain that
$\widetilde{f}_{D_i}(x_k^{p-1}x_hD_h)=0$ which completes the
inductive step. 
\qed

\subsection{Reduction to $W(n)_0$-invariant cohomology}

This subsection is devoted to prove the second step of the strategy
that was outlined in subsection 3.2. We consider the action of
$W(n)_{\geq 0}$ on $W(n)_{-1}$ obtained by the projection onto
$W(n)_0=W(n)_{\geq 0}/W(n)_{\geq 1}$ followed by the adjoint
representation of $W(n)_0=\mathfrak{gl}(n,F)$ on $W(n)_{-1}$.

\begin{pro}\label{W-rel-coho}
For every $s\in \Z_{\geq 0}$, we have
$$H^s(W(n),W(n)_{-1};W(n))=H^s(W(n)_{\geq 0}, W(n)_{-1}).$$
\end{pro}
\dem For every $s\in \Z_{\geq 0}$, consider the map
$$\phi_s: C^s(W(n),W(n)_{-1};W(n))
\to C^s(W(n)_{\geq 0},W(n)_{-1})$$ induced by the restriction to the
subalgebra $W(n)_{\geq 0}\subset W(n)$ and by the projection
$W(n)\twoheadrightarrow W(n)/W(n)_{\geq 0}=W(n)_{-1}$. It is
straightforward to check that the maps $\phi_s$ commute with the
differentials and hence they define a map of complexes. Moreover the
orthogonality conditions with respect to the subalgebra $W(n)_{-1}$
give the injectivity of the maps $\phi_s$. Indeed, on one hand, the
condition (\ref{ortho1}) says that an element $f\in
C^s(W(n),W(n)_{-1};W(n))$ is determined by its restriction to
$\wedge^s W(n)_{\geq 0}$. On the other hand, condition
(\ref{ortho2}) implies that the values of $f$ on an $s$-tuple are
determined, up to elements of $W(n)^{W(n)_{-1}}=W(n)_{-1}$, by
induction on the total degree of the $s$-tuple.

Therefore, to
conclude the proof, it is enough to prove that the maps $\phi_s$ are
surjective. Explicitly, if $f\in C^s(W(n)_{\geq 0}, W(n)_{-1})$,
consider the cochain $\widetilde{f}\in C^s(W(n), W(n))$ defined by
$$\widetilde{f}(x^{a^1},\cdots,x^{a^n})=\sum_{i=1}^n \sum_{\0\leq b^i <a^i}\prod_{i=1}^n
\binom{a^i}{b^i}f(x^{a^1-b^1},\cdots, x^{a^n-b^n})x^{b^1+\cdots +b^n} , $$
where if $a, b\in \N^n$ then $\binom{a}{b}:=\prod_{i=1}^n \binom{a_i}{b_i}$.

We are done if we show that $\widetilde{f} \in C^s(W(n),W(n)_{-1};W(n))$ since
it is clear that $\phi_s(\widetilde{f})=f$.  The first orthogonality condition
(\ref{ortho1}) follows easily from the definition. Consider the following expression
$$\widetilde{f}(x^{a^1},\cdots,D_j(x^{a^k}),\cdots,x^{a^n})=\sum_{\0\leq b'^k<a^k-\epsilon_j}\sum_{\stackrel{i\neq k}{\0\leq b^i< a^i}}
\left[(a^k)_j \binom{a^k-\ep_j}{b'^k}\prod_{i\neq k} \binom{a^i}{b^i}\right]$$
$$f(x^{a^1-b^1},\cdots, x^{a^k-\ep_j-b^k},\cdots, x^{a^n -b^n})x^{b^1+\cdots+
b'^k+\cdots+b^n}=$$
$$=\sum_{i=1}^n \sum_{\0\leq b^i <a^i}(b^k)_j\prod_{i=1}^n
\binom{a^i}{b^i}f(x^{a^1-b^1},\cdots, x^{a^n-b^n})x^{b^1+\cdots
+b^n-\ep_j}, $$ 
where we used the substitution $b^k=b'^k+\ep_j$ together with the equality
$(b^k)_j\binom{a^k}{b^k}=(a^k)_j\binom{a^k-\ep_j}{b^k-\ep_j}$.
Summing the above expression as $k$ varies from $1$ to $n$, we get
$[D_j,\widetilde{f}(x^{a^1},\cdots,x^{a^n})]$ which proves the
second orthogonality condition (\ref{ortho2}). \qed

\begin{pro}\label{W-H^1}
Consider $W(n)_{-1}$ as a trivial $W(n)_{\geq 1}$-module. Then
$$\begin{sis}
& H^1(W(n)_{\geq 0},W(n)_{-1})=0,\\
& H^2(W(n)_{\geq 0},W(n)_{-1})=H^2(W(n)_{\geq 1},W(n)_{-1})^{W(n)_0}.
\end{sis}$$
\end{pro}
\dem Consider the Hochschild-Serre spectral sequence relative to the
ideal \\ $W(n)_{\geq 1} \lhd W(n)_{\geq 0}$:
$$E_2^{r,s}=H^r(W(n)_0,H^s(W(n)_{\geq 1},W(n)_{-1}))\Rightarrow
H^{r+s}(W(n)_{\geq 0},W(n)_{-1}).$$ Note that since $T\subset
W(n)_{\geq 0}$, we can restrict to homogeneous cohomology (see
subsection $2.2$). Directly from homogeneity, it follows that the
first line $E_2^{*,0}=H^*(W(n)_0,W(n)_{-1})$ vanishes. Indeed the
weights that occur in $W(n)_{-1}$ are $-\epsilon_i$ while the
weights that occur in $W(n)_0$ are $\0$ and $\epsilon_i-\epsilon_j$.
Therefore the weights that occur in $W(n)_0^{\otimes k}$ have degree
congruent to $0$ modulo $p$ and hence they cannot be equal to
$-\epsilon_i$.

On the other hand, since $W(n)_{-1}$ is a trivial $W(n)_{\geq 1}$-module,
we have that
$$H^1(W(n)_{\geq 1}, W(n)_{-1})=\{f: W(n)_{\geq 1}\to W(n)_{-1} \: | \:
f([W(n)_{\geq 1}, W(n)_{\geq 1}])=0\}.$$
Therefore lemma \ref{W-commutators} gives that
$$H^1(W(n)_{\geq 1},W(n)_{-1})=\begin{sis}
&C^1(W(1)_1\oplus W(1)_2, W(n)_{-1}) \:\:\text{ if } n=1 \text{ and } p\geq 5,\\
&C^1(W(n)_1, W(n)_{-1})  \:\:\text{ if } n\geq 2 \text{ or } n=1 \text{ and }p=3.\\
\end{sis}$$
From this it follows that the second line $E_2^{*,1}=H^*(W(n)_{\geq
0},H^1(W(n)_{\geq 1}$, $W(n)_{-1}))$ vanishes again for homogeneity
reasons. Indeed, on one hand, the weights that appear on
$H^1(W(n)_{\geq 1},W(n)_{-1})$ have degree congruent to $2$ or $3$
modulo $p$ (the last one can occur only for $n=1$ and $p\geq 5$). On
the other hand the weights that appear on $W(n)_0$ (that are $\0$ or
$\epsilon_i-\epsilon_j$) are congruent to $0$ modulo $p$ and the
same is true for $W(n)^{\otimes k}$. \qed

\begin{lemma}\label{W-commutators}
Let $d\geq -1$ be an integer and suppose that it is different from $1$ if $n=1$.
Then $$[W(n)_1,W(n)_d]=W(n)_{d+1}.$$
\end{lemma}
\dem Clearly $[W(n)_1, W(n)_d]\subset W(n)_{d+1}$ by definition of
graded algebras. Consider formulas
$$ [x_i^2 D_i,x^bD_r]=\begin{cases}
b_i x^{b+\epsilon_i}D_r& \text{ if } i\neq r ,\\
(b_r-2) x^{b+\epsilon_r}D_r& \text{ if } i=r.
\end{cases} $$
Take an element $x^aD_r\in W(n)_{d+1}$. If $a_r\neq 0, 3$ the second formula above
with $i=r$ and $b=a-\epsilon_r$ shows that $x^aD_r \in [W(n)_1,W(n)_d]$.
On the other hand, if there exists some $i\neq r$ such that $a_i\neq  0,1$ then the first
formula above with $b=a-\epsilon_i$ gives that $x^aD_r \in [W(n)_1,W(n)_d]$.
Moreover if there is an index $s\neq r$ such that $a_s=1$, then we use the
formula
$$[x_s^2D_r, x^{a-\epsilon_s}D_s]=a_r x^{a-\epsilon_r+\epsilon_s}D_s-2 x^aD_r$$
since the first term on the right hand side
belongs to $[W(n)_1,W(n)_d]$ by what proved above.
Therefore, in virtue of our hypothesis on $d$, it remains to consider the elements $x_r^3D_r$
for $n\geq 2$. Choosing an $s\neq r$ we conclude by
$$[x_r^2D_s, x_rx_sD_r]=x_r^3 D_r-2 x_r^2 x_s D_s. \qed$$

\subsection{Computation of $W(n)_0$-invariant cohomology}

The aim of this subsection is to prove the following proposition that concludes
the third and last step of the proof.

\begin{pro}\label{W-H^2}
Denote with $\overline{{\rm Sq}(D_i)}$ the projection of ${\rm Sq}(D_i)$ onto $W(n)_{-1}$.
Then
$$H^2(W(n)_{\geq 1},W(n)_{-1})^{W(n)_0}=\bigoplus_{i=1}^n F \cdot \langle \overline
{{\rm Sq}(D_i)} \rangle,$$
with the exception of the case $n=1$ and $p=3$ when it is $0$.
\end{pro}

\dem First of all observe that if $n=1$ and $p=3$, then $W(n)_{\geq
1}=\langle x_1^2D_1 \rangle$ and hence the second cohomology group
vanishes. Hence we assume that $p\geq 5$ if $n=1$. It's easy to see
that the above cocycles $\overline{\Sq(D_i)}$ are $W(n)_0$-invariant
and independent modulo coboundaries (same argument as in section
$3.2$).
So we have to prove that they generate the second cohomology group. \\
Consider the truncated cohomology groups
$$H^2\left(\frac{W(n)_{\geq 1}}{W(n)_{\geq d}},W(n)_{-1}\right)^{W(n)_0}$$
as $d$ increases.
Observe that if $d\geq np-(n+1)$ then $W(n)_{\geq d+1}=0$ and
hence we get the cohomology we are interested in. Moreover if
$n\geq 2$ then lemma \ref{W-inv-(1,1)} below gives
$$H^2\left(\frac{W(n)_{\geq 1}}{W(n)_{\geq 2}},W(n)_{-1}\right)^{W(n)_0}=
C^2\left(W(n)_1,W(n)_{-1}\right)^{W(n)_0}=0, $$
while if $n=1$ (and $p\geq 5$) then by homogeneity we have that
$$H^2\left(\frac{W(1)_{\geq 1}}{W(1)_{\geq 3}},W(1)_{-1}\right)^{W(1)_0}=C^1\left(W(1)_1
\times W(1)_2,W(1)_{-1}\right)_{\0}=0.$$
The algebra $W(n)_{\geq 1}$ has a decreasing filtration
$\{W(n)_{\geq d}\}_{d=1,\cdots, n(p-1)-1}$ and the adjoint action of $W(n)_0$ respects
this filtration.
We consider one step of this filtration
\begin{equation*}
W(n)_d=\frac{W(n)_{\geq d}}{W(n)_{\geq d+1}}\lhd \frac{W(n)_{\geq 1}}{W(n)_{\geq d+1}}
\end{equation*}
and the related Hochschild-Serre spectral sequence
\begin{equation}\label{W-HocSer}
E_2^{r,s}=H^r\left(\frac{W(n)_{\geq 1}}{W(n)_{\geq d}},H^s\left(W(n)_d,W(n)_{-1}\right)\right)
\Rightarrow H^{r+s}\left(\frac{W(n)_{\geq 1}}{W(n)_{\geq d+1}},W(n)_{-1}\right).
\end{equation}
We fix a certain degree $d$ and we study, via the above spectral sequence, how the
truncated cohomology groups change if we pass from $d$ to $d+1$.
By what was said above, we can assume that $d>1$ if $n\geq 2$ and $d>2$ if $n=1$.

Observe that, since $W(n)_d$ is in the center of $W(n)_{\geq 1}/W(n)_{\geq d+1}$  and
$W(n)_{-1}$ is a trivial module, then $H^s\left(W(n)_d,W(n)_{-1}\right)=
C^s\left(W(n)_d,W(n)_{-1}\right)$ and  $W(n)_{\geq 1}/$ $ W(n)_{\geq d}$ acts trivially on it.
Since $E_{\infty}^{0,2}=0$ by lemma \ref{W-inv-(0,2)} below,
the above spectral sequence gives us the two following exact sequences
$$\xymatrix{
 C^1(W(n)_d,W(n)_{-1}) \ar@{^{(}->}[d]_{\alpha} & &  \\
 H^2\left(\frac{W(n)_{\geq 1}}{W(n)_{\geq d}},W(n)_{-1}\right) \ar@{->>}[d] & & \\
 E_{\infty}^{2,0} \ar@{^{(}->}[r] & H^2\left(\frac{W(n)_{\geq 1}}{W(n)_{\geq d+1}},
W(n)_{-1}\right) \ar@{->>}[r] & E_{\infty}^{1,1}    \\
}$$
where the injectivity of the map $\alpha$ follows from the exactness of the sequence
$$E_{\infty}^{1,0}=H^1\left(\frac{W(n)_{\geq 1}}{W(n)_{\geq d}},W(n)_{-1}\right)\hookrightarrow
 H^1\left(\frac{W(n)_{\geq 1}}
{W(n)_{\geq d+1}}, W(n)_{-1}\right)\twoheadrightarrow E_{\infty}^{0,1}={\rm Ker}(\alpha)$$
together with lemma \ref{W-commutators} which says that the first two terms are both equal to
$C^1\left(W(n)_1,W(n)_{-1}\right)$. Moreover, lemma \ref{W-commutators} gives that
\begin{equation}\label{W-(1,1)term}
E_{\infty}^{1,1}\subset E_2^{1,1}=
\begin{cases}
C^1(W(n)_1\times W(n)_d,W(n)_{-1}) & \text{ if } n\geq 2,\\
C^1\left([W(1)_1\oplus W(1)_2]\times W(1)_d,W(1)_{-1}\right) & \text{ if } n=1.
\end{cases}
\end{equation}
By taking cohomology with respect to $W(n)_0$ and using lemmas \ref{W-inv-(1,1)},
\ref{W-inv-(1,1)-n=1}, \ref{W-inv-(0,1)}, \ref{W-H^1-(0,1)} below,
 we see that the only terms responsible for the growth of the invariant
truncated cohomology groups are $H^1(W(n)_0,C^1(W(n)_d,W(n)_{-1}))$
if $n\geq 2$ and $d=p-1$ (see lemma \ref{W-H^1-(0,1)}) and
$(E_{\infty}^{1,1})^{W(n)_0}$ if $n=1$ and $d=p-2$ (see lemma
\ref{W-inv-(1,1)-n=1}). In both cases, we get the desired
statement. \qed

\begin{lemma}\label{W-inv-(0,2)}
In the above spectral sequence (\ref{W-HocSer}), we have $E_3^{0,2}=0$.
\end{lemma}
\dem By definition, $E_3^{0,2}$ is the kernel of the map
$$
\d: C^2\left(W(n)_d,W(n)_{-1}\right)=E_2^{0,2}\to E_2^{2,1}=
H^2\left(\frac{W(n)_{\geq 1}}{W(n)_{\geq d}},
C^1\left(W(n)_d,W(n)_{-1}\right)\right)
$$
that sends a $2$-cochain $f$ to the element $\d f$ given by
$\d f_{(E,F)}(G)=-f([E,F],G)$
whenever ${\rm deg}(E)+{\rm deg}(F)=d$ and $0$ otherwise.

The subspace of coboundaries $B^2\left(\frac{W(n)_{\geq 1}}{W(n)_{\geq d}},
C^1\left(W(n)_d, W(n)_{-1}\right)\right)$ is the image of the map
$$\partial: C^1\left(\frac{W(n)_{\geq 1}}{W(n)_{\geq d}}, C^1\left(W(n)_d,
W(n)_{-1}\right)\right)\to C^2\left(\frac{W(n)_{\geq 1}}{W(n)_{\geq d}}, C^1\left(W(n)_d,
W(n)_{-1}\right)\right)$$
that sends the element $g$ to the element $\partial g$ given by
$ \partial g_{(E,F)}(G)=-g_{[E,F]}(G).$
Hence $\partial g$ vanishes on the pairs $(E,F)$ for which ${\rm deg}(E)+{\rm deg}(F)=d$.

Therefore, if an element $f\in C^2\left(W(n)_d,W(n)_{-1}\right)$ is
in the kernel of $\d$, that is $\d f=\partial g$ for some $g$ as
before, then it should satisfy $f([E,F],G)=0$ for every $E, F, G$
such that ${\rm deg}(G)=d$ and ${\rm deg}(E)+{\rm deg}(F)=d$. By
letting $E$ vary in $W(n)_1$ and $F$ in $W(n)_{d-1}$, the bracket
$[E,F]$ varies in all $W(n)_d$ by lemma \ref{W-commutators} (note
that we are assuming $d\geq 3$ if $n=1$). Hence the preceding
condition implies that $f=0$. \qed

\begin{lemma}\label{W-inv-(1,1)}
If $n\geq 2$ and $d\geq 1$, then
$$C^1\left(W(n)_1\times W(n)_d,W(n)_{-1}\right)^{W(n)_0}=0.$$
\end{lemma}
\dem Note that invariance with respect to $T\subset W(n)_0$ is the
same as homogeneity, hence we can limit ourselves to considering
homogeneous cochains. In particular this implies the vanishing if
$d\not \equiv p-2 \mod p$.

Consider a homogeneous cochain $f\in C^1\left(W(n)_1\times W(n)_d,W(n)_{-1}\right)^{W(n)_0}$.
Since the action of $W(n)_0$ on $W(n)_1$ is transitive, the result will follow if we prove
that $f(x_1^2D_2,-)=0$. Indeed, assuming this is the case, imposing invariance respect
to an element $x_iD_j\in W(n)_0$, we get
$$0=(x_iD_j\circ f)(x_1^2D_2,-)=-f([x_iD_j,x_1^2D_2],-)-f(x_1^2D_2,[x_iD_j,-])+$$
$$+[x_iD_j, f(x_1^2D_2,-)]=-f([x_iD_j,x_1^2D_2],-),$$
which shows the vanishing for $f$ when restricted to
$[x_iD_j,x_1^2D_2]$. Continuing in this way one gets the vanishing
of $f$ on every element of $W(n)_1$ and hence the vanishing of $f$.
So it is enough to prove that for every element $x^aD_r\in W(n)_d$
one has $f(x_1^2D_2,x^aD_r)=0$.

Suppose that $p\geq 5$. Then by the homogeneity assumption on $f$,
we have the required vanishing as soon as $a_1=0$ or $a_2=p-1$
(because $p\geq 5!$). If $a_1\geq 1$ and $a_2<p-1$, we proceed by
induction on $a_1$. Suppose that we have proved the vanishing for
all the elements $x^bD_s$ such that $b_1<a_1$. Then, using the
induction hypothesis, the following invariance condition
$$0=(x_1D_2\circ f)(x_1^2D_2, x^{a-\epsilon_1+\epsilon_2}D_r)=
-f(x_1^2D_2,(a_2+1)x^aD_r)$$
gives the required vanishing.

Finally, in the case $p=3$, we can apply the same inductive argument, provided that we
first prove the vanishing in the case when $a_1=0$ or $a_2=p-1=2$.  This vanishing is provided by
the homogeneity of $f$ unless $x^aD_r$ is equal to $x_2^2D_2$, $x_2x_jD_j$ or
$x_1x_2^2x_j^2D_2$ (with $3\leq j\leq n$).
In this three exceptional cases one proves the vanishing using the
following invariance conditions:
$$\begin{sis}
&0=(x_1D_2\circ f)(x_1x_2D_2, x_2^2D_2)=-f(x_1^2D_2,x_2^2D_2)-f(x_1x_2D_2,2x_1x_2D_2),\\
&0=(x_jD_2\circ f)(x_1^2D_2,x_2^2D_j)=-f(x_1^2D_2,2x_2x_jD_j-x_2^2D_2),\\
&0=(x_jD_2\circ f)(x_1^2D_2,x_1x_2^2x_j^2D_2)=[x_jD_2, f(x_1^2D_2,x_1x_2^2x_j^2D_2)].
\end{sis} $$
\qed

\begin{lemma}\label{W-inv-(1,1)-n=1}
Consider the above spectral sequence (\ref{W-HocSer}). If $n=1$ then
$$(E_3^{1,1})_{\0}=(E_{\infty}^{1,1})_{\0}=\begin{cases}
\langle \overline{{\rm Sq}(D_1)} \rangle & \: \text{ if } d=p-2,\\
0& \: \text{ otherwise },
\end{cases}$$
where $\overline{{\rm Sq}(D_1)}$ denotes the restriction of
${\rm Sq}(D_1)$ to $W(1)_1\times W(1)_{p-2}$.
\end{lemma}
\dem For $n=1$ we have that $T=W(1)$ and therefore the
$W(n)_0$-invariance is the same as homogeneity. By Formula
(\ref{W-(1,1)term}) and homogeneity, we get
$$(E_2^{1,1})_{\0}=
\begin{cases}
\langle x_1^3D_1\times x_1^{p-2}D_1\to D_1\rangle & \: \text{ if } d=p-3 ,\\
\langle x_1^2D_1\times x_1^{p-1}D_1\to D_1\rangle & \: \text{ if } d=p-2,\\
0 & \: \text{ otherwise. }\\
\end{cases}$$
The term $(E_{\infty}^{1,1})_{\0}=(E_3^{1,1})_{\0}$ is the kernel of the differential map ${\rm d}:
(E_2^{1,1})_{\0}\to (E_2^{3,0})_{\0}=H^3\left(\frac{W(1)_{\geq 1}}
{W(1)_{\geq d}},W(1)_{-1}\right)_{\0}$. In view of the explicit description of $(E_2^{1,1})_{\0}$
as above, it's enough to show that the map ${\rm d}$ is different from $0$ if $d=p-3$, since
if $d=p-2$ then the cocycle $\overline{\Sq(D_1)}$ belongs to $(E_{\infty}^{1,1})_{\0}$ and is different from $0$ because ${\rm Sq}(D_1)(x_1^2D_1,x_1^{p-1}D_1)=-3 D_1$
($\neq 0$ for $p\geq 5$!).

So let $d=p-3$ (and hence $p\geq 7$) and suppose that
${\rm d}\langle x_1^3D_1\times x_1^{p-2}D_1\to D_1\rangle=\partial g$ for
$g\in C^2\left(\frac{W(1)_{\geq 1}}{W(1)_{\geq p-3}},W(1)_{-1}\right)_{\0}$.
If $p=7$ then $g=0$ for homogeneity reasons. Otherwise (if $p>7$) then note that
the cocycle ${\rm d}\langle x_1^3D_1\times x_1^{p-2}D_1\to D_1\rangle$ vanishes on the triples
$(x_1^2D_1,x_1^{j+1}D_1,x_1^{p-1-j}D_1)$ for $3\leq j \leq (p-3)/2$ and hence we get the following
 conditions on $g$:
\begin{equation*}
\begin{split}
0=\partial g(x_1^2D_1,x_1^{j+1}D_1,x_1^{p-1-j}D_1)=&-(j-1)g(x_1^{j+2}D_1,x_1^{p-1-j}D_1)+\\
&+(p-3-j)g(x_1^{p-j}D_1,x_1^{j+1}D_1)\\
\end{split}
\end{equation*}
from which, by decreasing induction on $j$, we deduce that $g(x_1^{j+1}D_1,x_1^{p-j}D_1)=0$
and hence that $g=0$. But this is absurd since
$${\rm d}\langle x_1^3D_1\times x_1^{p-2}D_1\to D_1\rangle(x_1^2D_1,x_1^3D_1,x_1^{p-3}D_1)=
-(p-5)D_1\neq 0. \qed $$

\begin{lemma}\label{W-inv-(0,1)}
Let $d\in \Z_{\geq 0}$. Then we have $C^1\left(W(n)_d,W(n)_{-1}\right)^{W(n)_0}=0$.
\end{lemma}
\dem Observe that $C^1\left(W(n)_d,W(n)_{-1}\right)^{W(n)_0}\subset
C^1\left(W(n)_d,W(n)_{-1}\right)_{\0}$ and the last term is non-
vanishing only if $d=p-1$ and $n\geq 2$, in which case we have the
homogeneous cochains $g(x_i^{p-1}x_jD_j)=a_j^i D_i$, $a_j^i\in F$
(for $i\neq j$). We get the vanishing of $g$ by means of the
following cocycle condition
\begin{equation}\label{W-coboundary}
0=\d g_{x_iD_j}(x_i^{p-2}x_j^2D_j)=-2g(x_i^{p-1}x_jD_j)=-2a_j^i D_i. \qed
\end{equation}

\begin{lemma}\label{W-H^1-(0,1)}
Let $d\in \Z_{\geq 0}$. Then
$$H^1\left(W(n)_0,C^1\left(W(n)_d,W(n)_{-1}\right)\right)=
\begin{cases}
\oplus_{i=1}^n \langle \overline{{\rm Sq}(D_i)} \rangle &  \text{if } n\geq 2
\, \text{ and } d=p-1,\\
0 &  \text{otherwise,}
\end{cases}$$
where $\overline{{\rm Sq}(D_i)}$ denotes the restriction of ${\rm Sq}(D_i)$ to
$W(n)_0\times W(n)_{p-1}$.
\end{lemma}
\dem Observe that, since the maximal torus $T$ is contained in
$W(n)_0$, the cohomology with respect to $W(n)_0$ reduces to
homogeneous cohomology. Hence the required group can be non-zero
only if $d\equiv p-1 \mod p$ (and hence only if $n\geq 2$). More
precisely, since the weights appearing on $W(n)_{-1}$ are
$-\epsilon_k$ and the weights appearing on $W(n)_0$ are
$\epsilon_i-\epsilon_j$ (possibly with $i=j$), the weights appearing
on $W(n)_d$ can be $-\epsilon_i+\epsilon_j -\epsilon_k$ (for every
$1\leq i, j, k\leq n$). Hence the required group can be non-zero
only if $d=p-1$ or $d=2p-1$ (this last case only if $n\geq 3$).

Consider first the case $d=2p-1$ ($n\geq 3$). A homogeneous cochain
$f\in C^1\left(W(n)_0,C^1\left(W(n)_{2p-1},W(n)_{-1}\right)\right)_{\0}$
takes the following non-zero values
$$f_{x_iD_j}(x_i^{p-1}x_k^{p-1}x_jx_hD_h)=\alpha_{ijk}^hD_k$$
for every $i,j,k$ mutually distinct and $h\neq i,k$. From the vanishing of $\d f$, we get
$$
0=\d f_{(x_iD_j,x_kD_i)}(x_i^{p-1}x_k^{p-1}x_jx_hD_h)=-[x_kD_i,f_{x_iD_j}(x_i^{p-1}
x_k^{p-1}x_jx_hD_h)]+$$
$$+f_{x_kD_j}(x_i^{p-1}x_k^{p-1}x_jx_hD_h)=\alpha_{ijk}^h D_i+\alpha_{kji}^hD_i ,
$$
$$
0=\d f_{(x_iD_jx_kD_j)}(x_i^{p-1}x_k^{p-1}x_jx_hD_h)=[x_iD_j,f_{x_kD_j}
(x_i^{p-1}x_k^{p-1}x_jx_hD_h)]
+$$
$$-[x_kD_j,f_{x_iD_j}(x_i^{p-1}x_k^{p-1}x_jx_hD_h)]=-\alpha_{kji}^hD_j+\alpha_{ijk}^hD_j.
$$
Adding these two equations, it follows that $2 \alpha_{ijk}^h=0$
and hence $f=0$.

Consider now the case $d=p-1$. First of all, a homogeneous cocycle $f$
must satisfy
$f_{x_iD_i}=0.$
Indeed, by formula (\ref{restriction}), we have $0=\d f_{| x_iD_i}=x_iD_i\circ f -
\d(f_{| x_iD_i})$ from which, since the first term vanishes for homogeneity reasons,
it follows that
$f_{| x_iD_i}\in C^1\left(W(n)_{p-1},W(n)_{-1}\right)^{W(n)_0}$ which is zero by
lemma \ref{W-inv-(0,1)}.  Therefore a homogeneous cocycle can take the following
non-zero values (for $i, j, k$ mutually distinct):
$$\begin{sis}
&f_{x_iD_j}(x_i^{p-2}x_j^2D_j)=\alpha_{ij} D_i, \\
&f_{x_iD_j}(x_i^{p-1}x_kD_k)=\alpha_{ij}^k D_j, \\
&f_{x_iD_j}(x_i^{p-2}x_jx_kD_k)=\beta_{ij}^k D_i, \\
&f_{x_iD_j}(x_k^{p-1}x_jD_i)=\gamma_{ij}^k D_k, \\
&f_{x_iD_j}(x_i^{p-1}x_jD_k)=\delta_{ij}^k D_k, \\
&f_{x_iD_j}(x_i^{p-1}x_jD_i)=\beta_{ij} D_i,  \\
&f_{x_iD_j}(x_i^{p-1}x_jD_j)=\gamma_{ij}D_j.
\end{sis}$$
By possibly modifying $f$ with a coboundary (see formula (\ref{W-coboundary})),
we can assume that $\underline{\alpha_{i,j}=0}$. Using this, we get the vanishing
of $\underline{\alpha_{ij}^k}$, $\underline{\beta_{ij}^k}$ and $\underline{\gamma_{ij}^k}$
by means of the following three cocycle conditions:
$$0=\d f_{(x_iD_j,x_iD_k)}(x_i^{p-2}x_k^2D_k)=[x_iD_j,f_{x_iD_k}(x_i^{p-2}x_k^2D_k)]
+f_{x_iD_j}(2x_i^{p-1}x_kD_k)=$$
$$=[-\alpha_{ik} +2 \alpha_{ij}^k] D_j,$$
$$0=\d f_{(x_iD_j,x_iD_k)}(x_i^{p-3}x_jx_k^2D_k)=-f_{x_iD_k}(x_i^{p-2}x_k^2D_k)
+f_{x_iD_j}(2x_i^{p-2}x_jx_kD_k)=$$
$$=[-\alpha_{ik}+2\beta_{ij}^k]D_i,$$
$$0=\d f_{(x_iD_j,x_kD_j)}(x_k^{p-2}x_j^2D_i)=-f_{x_kD_j}([x_iD_j,x_k^{p-2}x_j^2D_i])+$$
$$+f_{x_iD_j}([x_kD_j,x_k^{p-2}x_j^2D_i])=-f_{x_kD_j}(2x_k^{p-2}x_jx_iD_i)+f_{x_kD_j}(x_k^{p-2}x_j^2D_j)+$$
$$+f_{x_iD_j}(2x_k^{p-1}x_jD_i)=[-2\beta_{kj}^i+\alpha_{kj}+2\gamma_{ij}^k]D_k.$$
The coefficients $\underline{\delta_{ij}^k}$ and $\underline{\beta_{ij}}$ are
determined by the coefficients $\gamma_{ij}$ by the following two cocycle conditions:
$$0=\d f_{(x_iD_j,x_iD_k)}(x_i^{p-2}x_jx_kD_k)=-f_{x_iD_k}(x_i^{p-1}x_kD_k)+
f_{x_iD_j}(x_i^{p-1}x_jD_k)+$$
\begin{equation*}
-[x_iD_k,f_{x_iD_j}(x_i^{p-2}x_jx_kD_k)]=
[-\gamma_{ik}+\delta_{ij}^k+\beta_{ij}^k]D_k=
[-\gamma_{ik}+\delta_{ij}^k]D_k, \tag{*}
\end{equation*}
$$ 0=\d f_{(x_iD_j,x_jD_i)}(x_i^{p-1}x_jD_j)=f_{x_iD_j}(-x_i^{p-2}x_j^2D_j)+f_{x_iD_j}(-x_i^{p-1}x_jD_i)+$$
\begin{equation*}
-[x_jD_i,f_{x_iD_j}(x_i^{p-1}x_jD_j)]=[-\alpha_{ij}-\beta_{ij}+
\gamma_{ij}] D_i=[-\beta_{ij}+ \gamma_{ij}] D_i.  \tag{**}
\end{equation*}
The coefficients $\underline{\gamma_{ij}}$ satisfy the relation $\gamma_{ij}=
\gamma_{ik}$ (for $i, j, k$ mutually distinct as before). Indeed from the
cocycle condition
$$0=\d f_{(x_iD_j,x_iD_k)}(x_i^{p-1}x_jD_i)=-f_{x_iD_k}(x_i^{p-1}x_jD_j)
-[x_iD_k,f_{x_iD_j}(x_i^{p-1}x_jD_i)]+$$
$$+f_{x_iD_j}(-x_i^{p-1}x_jD_k)=[-\alpha_{ik}^j+\beta_{ij}-\delta_{ij}^k]D_k=
[\beta_{ij}-\delta_{ij}^k]D_k,$$
and using the relations $(*)$ and $(**)$ as above, we get
$\gamma_{ij}=\beta_{ij}=\delta_{ij}^k=\gamma_{ik}:=\gamma_i$.

We conclude the proof by observing that the elements $\overline{\Sq(D_i)}$ are independent
modulo coboundaries (if $n\geq 2$) as it follows from
$${\rm Sq}(D_i)(x_iD_j,x_i^{p-1}x_jD_j)=\frac{1}{(p-1)!}[D_i(x_iD_j),
(D_i)^{p-1}(x_i^{p-1}x_jD_j)]=D_j. $$ \qed

\section{The Special algebra}

\subsection{Definition and basic properties}

Throughout this section, we use the notations introduced in subsection $3.1$ and we
fix an integer $n\geq 3$.
Consider the following map, called divergence:

$${\rm div}:\begin{sis}
 W(n) & \to A(n)\\
 \sum_{i=1}^n f_iD_i & \mapsto \sum_{i=1}^n D_i(f_i).
\end{sis}$$
Clearly it is a linear map of degree $0$ that satisfies the following formula
(see \cite[chap. 4, lemma 3.1]{FS}):
\begin{equation*}
{\rm div}([D,E])=D({\rm div}(E))-E({\rm div}(D)).
\end{equation*}
Therefore the space $S'(n):=\{E\in W(n)\: | \: {\rm div}(E)=0\}$ is a graded subalgebra of $W(n)$
and we have an exact sequence of $S'(n)$-modules
\begin{equation}\label{S'-W-sequence}
0\to S'(n)\longrightarrow W(n)\stackrel{\rm div}{\longrightarrow} A(n)_{<\tau}\to 0.
\end{equation}

\begin{defi}
The Special algebra is the derived algebra of $S'(n)$:
$$S(n):=S'(n)^{(1)}=[S'(n),S'(n)].$$
\end{defi}

In order to describe the structure of $S(n)$, we introduce the following maps
(for $1\leq i,  j\leq n$)
$$D_{ij}:\begin{sis}
A(n)& \longrightarrow W(n)\\
f & \mapsto D_j(f)D_i-D_i(f)D_j.\\
\end{sis}$$
Note that $D_{ij}(A(n))\subset S'(n)$ and moreover if $\sum_{i=1}^n f_iD_i$ and
$\sum_{j=1}^ng_jD_j$ are two elements of $S'(n)$ then we have the following formula
\begin{equation}\label{S-bracket}
\left[\sum_{i=1}^n f_iD_i,\sum_{j=1}^ng_jD_j\right]=-\sum_{1\leq i, j\leq n}D_{ij}(f_ig_j),
\end{equation}
which in particular gives the following special case
\begin{equation}\label{D_ij-bracket}
[D_{ij}(f),D_{ij}(g)]=D_{ij}(D_{ij}(f)(g)).
\end{equation}

\begin{theo}\label{S-structure}
The algebra $S(n)$ satisfies the following properties:
\begin{itemize}
\item[(i)] $S(n)$ is generated by the elements $D_{ij}(f)$ for $f\in A(n)$
and $1\leq i<j\leq n$.
\item[(ii)] We have the following exact sequence of $S(n)$-modules
\begin{equation}\label{S-S'-sequence}
0\to S(n)\to S'(n)\to \oplus_{i=1}^n F\cdot \langle x^{\tau-(p-1)\epsilon_i}D_i
\rangle \to 0,
\end{equation}
where the last term is a trivial $S(n)$-module.
\item[(iii)] $S(n)$ is a restricted simple graded Lie algebra of dimension $(n-1)(p^n-1)$.
\end{itemize}
\end{theo}
\dem See \cite[Chap. 4, Prop. 3.3, Theo. 3.5 and 3.7]{FS}. \qed

Note that the unique term of negative degree is $S(n)_{-1}=\oplus_{i=1}^n
F\cdot \langle D_i \rangle$
while the term of degree $0$ is $S(n)_0=\oplus_{i=2}^n F \cdot
\langle x_iD_i-x_1D_1\rangle \oplus_{1\leq j\neq k\leq n}F\cdot \langle x_jD_k
\rangle$ and its adjoint action on $S(n)_{-1}$ induces
an isomorphism $S(n)_0\cong \mathfrak{sl}(n,F)$.

The algebra $S(n)$ admits a root space decomposition with respect to a canonical
Cartan subalgebra.

\begin{pro}\label{Cartan-S}
Recall that $h_i:=x_iD_i$ for every $i\in \{1,\cdots, n\}$.
\begin{itemize}
\item[(a)] $T_S:=T\cap S(n)=\oplus_{i=2}^n F\cdot \langle h_i-h_1 \rangle$
is a maximal torus of $H(n)$ (called the canonical maximal torus).
\item[(b)] The centralizer of $T_S$ inside $S(n)$ is the subalgebra
$$C_S=\bigoplus_{\stackrel{2\leq j\leq n}{0\leq a\leq p-2}}F \cdot\langle
D_{1j}(x^{\underline{a}+\epsilon_1+\epsilon_j})\rangle$$
which is hence a Cartan subalgebra (called the canonical Cartan subalgebra).
The dimension of $C_H$ is $(n-1)(p-1)$.
\item[(c)] Let $\Phi_S:={\rm Hom}_{\mathbb{F}_p}(\oplus_{i=2}^n \mathbb{F}_p
\langle h_i-h_1 \rangle , \mathbb{F}_p)$, where $\mathbb{F}_p$ is the prime
field of $F$. In the Cartan decomposition $S(n)=C_S\oplus_{\phi\in \Phi_S-\0}S(n)_{\phi}$,
the dimension of every $S(n)_{\phi}$, with $\phi\in \Phi_S-\0$, is $(n-1)p$.
\end{itemize}
\end{pro}
\dem See \cite[Chap. 4, Theo. 3.6]{FS}. \qed

\subsection{Strategy of the proof of the Main Theorem}

In this subsection, we outline the strategy of the proof of Theorem
\ref{S-finaltheorem} from the Introduction. Hence, from now on, we assume that
the characteristic $p$ of the base field $F$ is different from $2$.

We first check that $\Theta$ is a cocycle. It is enough to verify that it is a cocycle when
restricted to $S(n)_{-1}$ and that it is $S(n)_0$-invariant:
$$\d \Theta(D_i,D_j,D_k)=[D_i,D_{jk}(x^{\tau})]-[D_j,D_{ik}(x^{\tau})]+[D_k,D_{ij}(x^{\tau})]=$$
$$= -D_{jk}(x^{\tau-\epsilon_i})+D_{ik}(x^{\tau-\epsilon_j})-D_{ij}(x^{\tau-\epsilon_k})=0$$
and (for $h\neq k$)
$$(x_hD_k\circ \Theta)(D_i,D_j)=[x_hD_k,D_{ij}(x^{\tau})]+\delta_{ih}\Theta(D_k,D_j)
+\delta_{jh}\Theta(D_i,D_k) =$$
$$=\delta_{hj}D_{ki}(x^{\tau})-\delta_{hi}D_{kj}(x^{\tau})+\delta_{ih}D_{kj}(x^{\tau})
+\delta_{jh}D_{ik}(x^{\tau})=0.$$

Moreover the cocycles $\Theta$ and $\Sq(D_i)$ appearing in Theorem \ref{S-finaltheorem}
are independent modulo coboundaries. Indeed, if
$\gamma\in \{{\rm Sq}(D_1),\cdots,{\rm Sq}(D_n),\Theta\}$ then we have (for $i\neq j$)
$$\gamma(D_i,D_j)=
\begin{cases}
D_{ij}(x^{\tau}) &\text{ if } \gamma=\Theta, \\
0 &\text{ otherwise, }
\end{cases} \hspace{1cm} \text{ and }$$
\begin{equation}\label{S-ind-coc}
\gamma(x_iD_j,D_{ji}(x_i^{p-1}x_j^2))=
\begin{cases}
-2 D_i & \text{ if } \gamma={\rm Sq}(D_i) , \\
0  & \text{ otherwise,}
\end{cases}
\end{equation}
while for every $g\in C^1(S(n),S(n))$ the coboundary $\d g(D_i,D_j)=[D_i,g(D_j)]-[D_j,g(D_i)]$
cannot contain the monomial $D_{ij}(x^{\tau})$ for degree reasons and \\
$\d g(x_iD_j, D_{ji}(x_i^{p-1}x_j^2))=[x_iD_j,g(D_{ji}(x_i^{p-1}x_j^2))]-[D_{ji}(x_i^{p-1}x_j^2),g(x_iD_j)]$
cannot contain the monomial $D_i$.

Assuming the results of the next subsection, we complete the proof of the
Theorem \ref{S-finaltheorem}.

\dem[Proof of Theorem \ref{S-finaltheorem}] From the sequence
(\ref{S'-W-sequence}), using proposition \ref{S-H^1}, we get the
exact sequence
$$0 \to H^1(S(n),A(n)_{<\tau}) \stackrel{\partial}{\longrightarrow} H^2(S(n),S'(n))\to
H^2(S(n),W(n)).$$
By proposition \ref{S-H^2}, we known that $H^2(S(n),W(n))$ is generated by
the cocycles ${\rm Sq}(D_i)$. These clearly belong to $H^2(S(n),S'(n))$ and hence
the above exact sequence splits
$$H^2(S(n),S'(n))=  \bigoplus_{i=1}^n  \langle {\rm Sq}(D_i)\rangle \bigoplus \partial
H^1(S(n),A(n)_{<\tau}).$$
On the other hand, from the sequence (\ref{S-S'-sequence}), we get the
exact sequence
$$0\to H^2(S(n),S(n))\to H^2(S(n),S'(n)) \to \bigoplus_{i=1}^n
H^2(S(n),x^{\tau-(p-1)\epsilon_i}D_i),$$
where we used that $H^1(S(n),M)=0$ for a trivial $S(n)$-module $M$.
Since the cocycles ${\rm Sq}(D_i)$ belong to $H^2(S(n),S(n))$, we are left with verifying
which of the elements of $\partial H^1(S(n),A(n)_{<\tau})$ (which we know by proposition
\ref{S-coho-A}) belong to $H^2(S(n),S(n))$.

Consider first the cocycle ${\rm ad}(x^{\tau}): D_i\mapsto D_i(x^{\tau})=-x^{\tau-\epsilon_i}$.
It lifts to the cocycle $\widetilde{{\rm ad}}(x^{\tau})\in C^1(S(n),W(n))$ given by
$\widetilde{{\rm ad}}(x^{\tau}):D_i\mapsto x^{\tau}D_i$ and $0$ on the other elements.
Therefore the only non-zero values of $\partial (\widetilde{{\rm ad}}(x^{\tau}))$ can be
(for $k\neq h$):
$$\partial (\widetilde{{\rm ad}}(x^{\tau}))(D_i,D_j)=[D_i,x^{\tau}D_j]-[D_j,x^{\tau}D_i]=-D_{ij}
(x^{\tau})$$
$$\partial (\widetilde{{\rm ad}}(x^{\tau}))(D_i,x_kD_h)=-[x_kD_h,x^{\tau}D_i]-
\widetilde{{\rm ad}}(x^{\tau})([D_i,x_kD_h])=$$
$$=\delta_{ik}x^{\tau}D_h-\phi(\delta_{ik}D_h)=0$$
and hence we have that $\partial({\rm ad}(x^{\tau}))=-\Theta$.

Consider now the element $\chi_i\in H^1(S(n),A(n)_{<\tau})$ and
choose a lifting $\widetilde{\chi_i}\in$ $ C^1(S(n),W(n))$ in such a
way that if $\chi_i(\gamma)=0$ then $\widetilde{\chi_i}(\gamma)=0$.
Then (if $j\neq i$), we have
$$\partial (\chi_i)(D_j,x^{\tau-(p-1)(\epsilon_i+\epsilon_j)}D_i)=
[D_j,\widetilde{\chi_i}(x^{\tau-(p-1)(\epsilon_j+\epsilon_i)})]=x^{\tau-(p-1)\epsilon_j}D_j,$$
because the only possible lifting to $W(n)$ of the element
$\chi_i(x^{\tau-(p-1)(\epsilon_i+\epsilon_j)}D_i)=x^{\tau-(p-1)\epsilon_j}$
is $x^{\tau-(p-2)\epsilon_j}D_j$. On the other hand, for every
cochain $g\in C^1(S(n),$$ \\$$x^{\tau-(p-1)\epsilon_j})$ we have $\d
g(D_j,x^{\tau-(p-1)\epsilon_i-(p-1)\epsilon_j}D_i)=0$ because the
module is trivial and
$[D_j,x^{\tau-(p-1)\epsilon_i-(p-1)\epsilon_j}D_i]=0$. Hence the
projection of $\partial(\chi_i)$ into
$H^2(S(n),$ $ x^{\tau-(p-1)\epsilon_j})$ is non-zero and therefore
$\partial(\chi_i)\not\in H^2(S(n),S(n))$. \qed

\begin{pro}\label{S-coho-A}
Consider the natural action of $S(n)$ on $A(n)_{<\tau}$. We have
$$H^1(S(n),A(n)_{<\tau})=\oplus_{i=1}^n \langle \chi_i \rangle \oplus \langle {\rm ad}(x^{\tau})
\rangle, $$
where the $\chi_i\in H^1(S(n),A(n)_{<\tau})$ are defined by
$$\chi_i(x^aD_k)=\begin{cases}
x^a\cdot x_i^{p-1} & \text{ if } k=i,\\
0 & \text{ otherwise. }
\end{cases}$$
\end{pro}
\dem First of all note that $\chi_i$ takes values in $A(n)_{<\tau}$
(and not merely on $A(n)$) since $x^{\tau-(p-1)\epsilon_i}D_i\not\in
S(n)$. To prove that $\chi_i$ are cocycles, it is enough to verify
the following two cocycle conditions (where $j, h, k$ are different
from $i$)
$$
{\rm d}(\chi_i)(D_{ij}(x^a),D_{hk}(x^b))=-D_{hk}(x^b)(D_j(x^a)x_i^{p-1})-
\chi_i([D_{ij}(x^a),D_{hk}(x^b)])=$$
$$=-D_{hk}(x^b)(D_j(x^a)x_i^{p-1}+\chi_i(D_{hk}(x^b)(D_j(x^a))D_i)=0,$$
$${\rm d}(\chi_i)(D_{ij}(x^a),D_{ih}(x^b))=D_{ij}(x^a)(D_h(x^b)x_i^{p-1})-
D_{ih}(x^b)(D_j(x^a)x_i^{p-1})+$$
$$-\chi_i(D_{ij}(x^a)(D_h(x^b))D_i-D_{ih}(x^b)(D_j(x^a))D_i)=$$
$$=D_h(x^b)D_{ij}(x^a)(x_i^{p-1})-D_j(x^a)D_{ih}(x^b)(x_i^{p-1})=0.$$

The independence of the above cocycles $\gamma_i$ and ${\rm ad}(x^{\tau})$
modulo coboundaries follows from the fact that if
$\gamma\in\{\chi_1,\cdots,\chi_n,{\rm ad}(x^{(\tau)})\}$ then
$$\gamma(D_i)=\begin{cases}
x_i^{p-1} & \text{ if } \gamma=\chi_i, \\
D_i(x^{\tau})=-x^{\tau-\epsilon_i} & \text{ if } \gamma={\rm ad}(x^{\tau}), \\
0 & \text{ otherwise, }
\end{cases}$$
while for any $g\in A(n)_{<\tau}$ the element ${\rm d}g(D_i)=D_i(g)$ cannot the
monomials $x_i^{p-1}$
or $x^{\tau-\epsilon_i}$.

In order to prove that the whole cohomology group is generated by the above
cocycles, we consider the exact sequence of $S(n)$-modules
$$0\to A(n)_{<\tau}\to A(n)\to \langle x^{\tau} \rangle \to 0 ,$$
where $\langle x^{\tau}\rangle$ is a trivial $S(n)$-module. By taking cohomology and using
the fact that $H^1(S(n),x^{\tau})=0$, we obtain
$$H^1(S(n),A(n)_{<\tau})=\langle {\rm ad}(x^{\tau})\rangle \oplus H^1(S(n),A(n)).$$

Finally, to compute the last cohomology group we use the Hochschild-Serre spectral sequence
with respect to the subalgebra $S(n)_{-1}<S(n)$:
$$E_1^{r,s}=H^s(S(n)_{-1},C^r(S(n)/S(n)_{-1},A(n)))\Rightarrow H^{r+s}(S(n),A(n)).$$
Note that $E_1^{0,1}=H^1(S(n)_{-1},A(n))=\oplus_{i=1}^n F\cdot
\langle x_i^{p-1} D_i^*\rangle$ (by proposition
\ref{W(n)_{-1}-cohomology}) and the $\chi_i$ are global cocycles
lifting them. On the other hand, by the same argument as in
proposition \ref{W-rel-coho}, we have that
$E_2^{1,0}=H^1(S(n),S(n)_{-1};A(n))=H^1(S(n)_{\geq 0}, 1)$. But this
last group vanishes since $[S(n)_{\geq 0},S(n)_{\geq 0}]=S(n)_{\geq
0}$ as it follows easily from lemma \ref{S-commutators} above.
\qed

In the course of the proof of the main result, we obtain a new
proof of the following result.

\begin{theo}[Celousov]\label{S-Celousov}
$$H^1(S(n),S(n))=\oplus_{i=1}^n
{\rm ad}(x^{\tau-(p-1)\epsilon_i}D_i) \oplus {\rm ad}(x_1D_1).$$
\end{theo}
\dem From the exact sequence (\ref{S-S'-sequence}) of $S(n)$-modules
and using the fact that
$S'(n)^{S(n)}=H^1(S(n),x^{\tau-(p-1)\epsilon_i})=0$, we get that
$$H^1(S(n),S(n))= \oplus_{i=1}^n \langle {\rm ad}(x^{\tau-(p-1)\epsilon_i}D_i)\rangle \oplus
H^1(S(n),S'(n)).$$ From the exact sequence (\ref{S'-W-sequence}) and
using the facts that $W(n)^{S(n)}=0$ and $A(n)_{<\tau}^{S(n)}=F
\cdot \langle 1 \rangle$ together with proposition \ref{S-H^1}, an
easy computation with the coboundary map gives
$H^1(S(n),S'(n))=\langle {\rm ad}(x_1D_1) \rangle.$ \qed

\subsection{Cohomology of $W(n)$}

In this section we complete the proof of the main Theorem
by computing the first and the second cohomology group of $W(n)$ as a
$S(n)$-module.

\begin{pro}\label{S-H^1}
$H^1(S(n),W(n))=0.$
\end{pro}
\dem Consider the homogeneous Hochschild-Serre spectral sequence
(\ref{HS-ss1}) with respect to the subalgebra $S(n)_{-1}< S(n)$:
\begin{equation}\label{S-HS1}
(E_1^{r,s})_{\0}=H^s\left(S(n)_{-1}, C^r\left(S(n)/S(n)_{-1},W(n)\right)\right)_{\0}
\Rightarrow H^{r+s}\left(S(n),W(n)\right)_{\0}.
\end{equation}
Note that the vertical line $E_1^{0,*}=H^*(S(n)_{-1},W(n))_{\0}=H^*(W(n)_{-1},W(n))_{\0}$ vanishes
by corollary \ref{W(n)_{-1}-vanishing} and hence we get that
$$H^1(S(n);W(n))=H^1(S(n),S(n)_{-1};W(n)).$$
The same argument of proposition \ref{W-rel-coho}, using $S(n)^{S(n)_{-1}}=S(n)_{-1}$, gives
that
$$H^1(S(n),S(n)_{-1};W(n))=H^1(S(n)_{\geq 0},S(n)_{-1}),$$
 where $S(n)_{-1}$ is a $S(n)_{\geq 0}$-module via the
projection $S(n)_{\geq 0}\twoheadrightarrow  S(n)_0$ followed
by the adjoint representation of $S(n)_0$ on $S(n)_{-1}$.

Now consider  the Hochschild-Serre spectral sequence (\ref{HS-ss2})
relative to the ideal $S(n)_{\geq 1}\lhd S(n)_{\geq 0}$:
\begin{equation}\label{S-HS2}
E_2^{r,s}=H^r(S(n)_0,H^s(S(n)_{\geq 1},S(n)_{-1}))\Rightarrow H^{r+s}(S(n)_{\geq 0},S(n)_{-1}).
\end{equation}
By direct inspection, it is easy to see that
$E_2^{1,0}=H^1(S(n)_0,S(n)_{-1})=0$ for homogeneity reasons. On the
other hand, since $S(n)_{-1}$ is a trivial $S(n)_{\geq 1}$-module,
it follows from lemma \ref{S-commutators} above that
$H^1(S(n)_{\geq 1},S(n)_{-1})=C^1(S(n)_1, S(n)_{-1})$ and hence that
$E_2^{0,1}=C^1(S(n)_1,S(n)_{-1})^{S(n)_0}=0$ by lemma
\ref{S-inv-(0,1)} above. \qed

\begin{lemma}\label{S-commutators}
Let $d\geq -1$ be an integer.
Then $$[S(n)_1,S(n)_d]=S(n)_{d+1}.$$
\end{lemma}
\dem The inclusion $[S(n)_1,S(n)_d]\subset S(n)_{d+1}$ is obvious,
so we fix an element $D_{ij}(x^a)\in S(n)_{d+1}$ (that is ${\rm
deg}(x^a)=d+3\geq 2$) and we want to prove that it belongs to
$[S(n)_1,S(n)_d]$.

Suppose first that $a_i\geq 2$ and $a_j<p-1$. Then we are done by formula
$$[x_i^2D_j,D_{ij}(x^{a-2\epsilon_i+\epsilon_j})]=
D_{ij}(x_i^2D_j(x^{a-2\epsilon_i+\epsilon_j}))=(a_j+1)D_{ij}(x^a). $$
Therefore (by interchanging $i$ and $j$) it remains to consider the elements $x^a$ for which
$a_i=a_j=p-1$ or $0\leq a_i,a_j\leq 1$. We first consider the elements satisfying this
latter possibility. If $a_i=a_j=1$ then we use formula (see
(\ref{D_ij-bracket}))
$$[D_{ij}(x_i^2x_j), D_{ij}(x^{a-\epsilon_i})]=D_{ij}((x_i^2D_i-2x_ix_jD_j)(x^{a-\epsilon_i}))=
-2D_{ij}(x^a).$$
On the other hand, if $(a_i,a_j)=(1,0)$ then, by the hypothesis ${\rm deg}(x^a)=d+3\geq 2$,
there should
exist an index $k\neq i,j$ such that $a_k\geq 1$ and hence we use formula
$$[D_{ij}(x_i^2x_k),D_{ij}(x^{a-\epsilon_i-\epsilon_k+\epsilon_j})]=-2D_{ij}(x^a).$$
Analogously, if $a_i=a_j=0$ then there should exist either two
different indices $k, h\not\in\{i,j\}$ such that $a_k,a_h\geq 1$
or one index $k\neq i,j$ such that $a_k\geq 2$. We reach the
desired conclusion using formula (with $h=k$ in the second case)
$$[D_{ij}(x_kx_hx_j), D_{ij}(x^{a-\epsilon_h-\epsilon_k+\epsilon_i})]=D_{ij}(x^a).$$
Hence we are reduced to considering the elements $D_{ij}(x^a)$ such that $a_i=a_j=p-1$.
Here we have to use the hypothesis that $n\geq 3$.
Suppose first that there exist an index $k\not\in\{i,j\}$ such that $a_k\neq p-2$.
Consider formula (see (\ref{S-bracket}))
$$[D_{ik}(x_kx_i^2),D_{ij}(x^{a-\epsilon_i})]=-2D_{ij}(x^a)+2 D_{ik}(x^{a+\epsilon_k-
\epsilon_j})+4 D_{kj}(x^{a-\epsilon_i+\epsilon_k}).$$
The last two elements have $k$-coefficients  different from $p-1$ (by the hypothesis
$a_k\neq p-2$) and therefore belong to $[S(n)_1,S(n)_d]$ by what proved above.
This implies also that our element $D_{ij}(x^a)$ belongs to $[S(n)_1,S(n)_d]$.

At this point, only the elements $D_{ij}(x^{a})$ with $a=\underline{p-2}+\ep_i+\ep_j$ are left.
Consider the following linear system (where $k\neq i, j$):
$$\begin{sis}
& [D_{ik}(x_kx_i^2), D_{ij}(x^{a-\ep_i})]=
-2D_{ij}(x^{a})+2D_{ik}(x^{a-\ep_j+\ep_k})-4D_{jk}(x^{a-\ep_i+\ep_k}),\\
& [D_{ik}(x_k^2x_i), D_{ij}(x^{a-\epsilon_k})]=-2D_{ij}(x^a)+
D_{ik}(x^{a-\epsilon_j+\epsilon_k})-D_{jk}(x^{a-\epsilon_i+\epsilon_k}),\\
& [D_{ik}(x_ix_jx_k), D_{ij}(x^{a-\epsilon_j})]=-D_{ij}(x^a)+
2D_{ik}(x^{a-\epsilon_j+\epsilon_k})-D_{jk}(x^{a-\epsilon_i+\epsilon_k}).
\end{sis}$$
Since the matrix
$\begin{pmatrix}
-2& 2 & -4\\
-2 & 1 & -1 \\
-1 & 2 & -1 \\
\end{pmatrix}$
has determinant equal to $8$ and hence is invertible over $F$, from
the preceding system we get that $D_{ij}(x^a)\in [S(n)_1,S(n)_d]$.
\qed

\begin{pro}\label{S-H^2}
Assume that the characteristic of the base field $F$ is different from $3$ if
$n=3$. Then
$$H^2(S(n),W(n))=\oplus_{i=1}^nF\cdot \langle {\rm Sq}(D_i)\rangle .$$
\end{pro}
\dem We have already proved that the above cocycles are independent
modulo coboundaries so that we are left with showing that they
generate the whole second cohomology group. This will be done in
several steps.\\
\underline{STEP I} : $H^2(S(n),W(n))=H^2(S(n),S(n)_{-1};W(n)).$

Consider the homogeneous Hochschild-Serre spectral sequence
(\ref{S-HS1}) with respect to the subalgebra $S(n)_{-1}< S(n)$.
Since, by corollary \ref{W(n)_{-1}-vanishing}, the vertical line
$E_1^{0,*}=$ $H^*(S(n)_{-1},$ $W(n))_{\0}$ vanishes, we will
conclude this first step by showing that $(E_2^{1,1})_{\0}=0$.

The proof of that is similar to the one of proposition \ref{WE_2^{1,1}}. We sketch a
proof referring to that proposition for notations and details. So suppose that we have
an element $[f]\in (E_1^{1,1})_{\0}=H^1(S(n)_{-1},C^1(S(n)/S(n)_{-1},W(n))_{\0}$
that goes to $0$ under the differential map ${\rm d}:(E_1^{1,1})_{\0}\to (E_1^{2,1})_{\0}$.
First of all, arguing by induction on degree as in proposition \ref{WE_2^{1,1}}, we can find a
representative $\widetilde{f}$ of the class $[f]$ such that for a certain $d$ and for
every $i=1,\cdots, n$, we have that
$$\begin{sis}
&\widetilde{f}_{D_i}(F)=0 & \quad \text{ for every } F \in S(n) \: : {\rm deg}(F)<d,\\
&\widetilde{f}_{D_i}(E)\in\langle x_i^{p-1}\rangle \otimes W(n)_{-1} & \quad \text{ for
every } E \in S(n) \: : {\rm deg}(E)=d.\\
\end{sis}$$
By homogeneity, it is easy to see that $\widetilde{f}_{D_i}$ can take non-zero values only on the
elements $E$ of the form (for a certain $k$)
$$\begin{sis}
&D_{kh}(x^{\underline{a}+\epsilon_h}) \text{ for } 1\leq a\leq p-1 \text{ and } h\neq k,
&\text{ (I) } \\
&x_k^{p-1} (x_rD_r-x_sD_s) \text{ for } k, r, s \text{ mutually distinct.}  &\hspace{1,6cm}
\text{ (II) } \\
\end{sis}$$
In particular, note that the degree $d$ of $E$ is at least $n-1\geq 2$. Now we can conclude the
proof using exactly
the same argument as in proposition \ref{WE_2^{1,1}}: we have to find,
for every $E$ as above, two elements
$A\in S(n)_0$ and $B\in S(n)_d$ such that $[A,B]=E$ and $A\not \in S(n)_{-\epsilon_j},
S(n)_{\epsilon_2+\cdots+\epsilon_n}$
for any $j=2,\cdots,n$ (which are exactly the weights appearing on $S(n)_{-1}$).
Explicitly: if $E$ is of type $(II)$ we take $B=x_kD_r$ and $A=1/2 \cdot
D_{rs}(x_k^{p-2}x_r^2x_s)$; if $E$ is of type
$(I)$ with $a\neq p-2$ then we take $B=x_kD_h$ and $A=-1/(a+2)\cdot D_{kh}(x^{\underline{a}+
2\epsilon_h-\epsilon_k})$.
Finally if $E$ is of type $(I)$ with $a=p-2$, then, choosing an index $j$ different from $k$ and
$h$ (this is possible since $n\geq 3$), the same argument as above gives the
vanishing of $\widetilde{f}_{D_i}$ on the following two elements
$$\begin{sis}
&3D_{hj}(x^{\underline{p-2}-\epsilon_k+\epsilon_j+\epsilon_h})-D_{hk}(x^{\underline{p-2}+
\epsilon_h})=
[x_kD_h,D_{jk}(x^{\underline{p-2}-\epsilon_k+\epsilon_j+\epsilon_h})],\\
&2D_{hj}(x^{\underline{p-2}-\epsilon_k+\epsilon_j+\epsilon_h})-2 D_{hk}(x^{\underline{p-2}+
\epsilon_h})=[x_jD_h,D_{jk}(x^{\underline{p-2}+\epsilon_h})].\\
\end{sis}$$
But then, since the matrix
$\left(\begin{array}{cc}
3& -1 \\
2 & -2 \\
\end{array}\right)$
has determinant equal to $-4$ and hence is invertible over $F$, we can take an appropriate linear
combination of the two elements above to get the vanishing of $\widetilde{f}_{D_i}$
on the element $D_{hk}(x^{\underline{p-2}+\epsilon_h})$.\\
\underline{STEP II} : $H^2(S(n),S(n)_{-1};W(n))\hookrightarrow
H^2(S(n)_{\geq 1},S(n)_{-1})^{S(n)_0}.$

First of all, exactly as in proposition \ref{W-rel-coho} (using that
$S(n)^{S(n)_{-1}}=S(n)_{-1}$), we get
$$H^2(S(n),S(n)_{-1};W(n))= H^2(S(n)_{\geq 0},S(n)_{-1})$$
where as usual $S(n)_{-1}$ is a $S(n)_{\geq 0}$-module via the
projection $S(n)_{\geq 0}\twoheadrightarrow  S(n)_0$ followed
by the adjoint representation of $S(n)_0$ on $S(n)_{-1}$.

Finally, we consider the Hochschild-Serre spectral sequence (\ref{S-HS2})
with respect to the ideal
$S(n)_{\geq 1}\lhd S(n)_{\geq 0}$.
Using that
$E_2^{2,0}=H^2(S(n)_0,S(n)_{-1})=0$ for homogeneity reasons and
$E_2^{1,1}=H^1(S(n)_0,C^1(S(n)_1,S(n)_{-1}))=0$ by lemmas \ref{S-commutators}
 and \ref{S-H^1-(0,1)}, we get the inclusion
$$H^2(S(n)_{\geq 0},S(n)_{-1})\hookrightarrow H^2(S(n)_{\geq 1},S(n)_{-1})^{S(n)_0}.$$
\underline{STEP III} : $ H^2(S(n)_{\geq 1},S(n)_{-1})^{S(n)_0}=\oplus_{i=1}^n
F\cdot \langle {\rm Sq}(D_i)\rangle.$

The strategy of the proof is the same as that of proposition \ref{W-H^2}: to compute,
step by step as $d$ increases, the truncated invariant cohomology groups
$$H^2\left(\frac{S(n)_{\geq 1}}{S(n)_{\geq d+1}},S(n)_{-1}\right)^{S(n)_0}.$$
By lemma \ref{S-inv-(1,1)}, we get that
$H^2\left(\frac{S(n)_{\geq 1}}{S(n)_{\geq 2}},1\right)^{S(n)_0}=
C^2\left(S(n)_1,S(n)_{-1}\right)^{K(n)_0}=0. $
On the other hand, if $d\geq n(p-1)-2$ then $S(n)_{\geq d+1}=0$ and
hence we get the cohomology we are interested in.

Consider the Hochschild-Serre spectral sequence associated to the ideal
$
S(n)_d=\frac{S(n)_{\geq d}}{S(n)_{\geq d+1}}\lhd \frac{S(n)_{\geq 1}}{S(n)_{\geq d+1}}
$:
\begin{equation}\label{S-HocSer}
E_2^{r,s}=H^r\left(\frac{S(n)_{\geq 1}}{S(n)_{\geq d}},H^s(S(n)_d,S(n)_{-1})\right)\Rightarrow
H^{r+s}\left(\frac{S(n)_{\geq 1}}{S(n)_{\geq d+1}},S(n)_{-1}\right).
\end{equation}
We get the same diagram as in proposition \ref{W-H^2} (the vanishing
of $E_3^{0,2}$ and the injectivity of the map $\alpha$ are proved in
exactly the same way). We conclude by taking cohomology with respect
to $S(n)_0$ and using lemmas \ref{S-inv-(1,1)},
\ref{S-inv-(0,1)} and \ref{S-H^1-(0,1)} below. \qed

\begin{lemma}\label{S-inv-(1,1)}
Assume that the characteristic of $F$ is different from $3$ if $n=3$. Then in the above
spectral sequence (\ref{S-HocSer}), we have that
$$(E_{\infty}^{1,1})^{S(n)_0}=0.$$
\end{lemma}
\dem For the above spectral sequence (\ref{S-HocSer}), we have the
inclusion
$$(E_{\infty}^{1,1})^{S(n)_0}\subset (E_2^{1,1})^{S(n)_0}=C^1(S(n)_1\times S(n)_d,S(n)_{-1})^{S(n)_0}.$$
Let $f$ be a homogeneous cochain belonging to $C^1(S(n)_1\times S(n)_d,S(n)_{-1})^{S(n)_0}$.
Since the action of $S(n)_0$ on S$(n)_1$ is transitive, the cochain $f$ is determined by
its restriction
$f(x_1^2D_2,-)$
(see the proof of lemma \ref{W-inv-(1,1)}). Even more, $f$ is determined by its restriction
to the pairs $(x_1^2 D_2, E)$ for which
$f(x_1^2D_2,E)\in \langle D_2\rangle$, which is equivalent to
$E\in S(n)_{-2\sum_{i\geq 2}
\epsilon_i}$ by the homogeneity of $f$.
Indeed, the values of $f$ on the other pairs $(x_1^2, F)$ for which
$f(x_1^2D_2,F)\in \langle D_j\rangle$ (for a certain $j\neq 2$) are determined by
the invariance condition
$$0=(x_jD_2\circ f)(x_1^2D_2, F)=[x_jD_2,f(x_1^2D_2,F)]-f(x_1^2D_2,[x_jD_2,F]).$$
A base for the space $S(n)_{-2\sum_{i\geq 2}\epsilon_i}$ consists of the elements
\begin{align*}
& D_{1k}(x^{\underline{p-1}-\epsilon_1+\epsilon_k}) \: \text{ for }\: k\neq 1, \tag{A}\\
& D_{3h}(x^{\underline{a}-2\epsilon_1+\epsilon_3+\epsilon_h}) \:\text{ for }\: 0\leq a\leq p-2
\: \text{ and } h\neq 3. \tag{B}
\end{align*}
For the elements of type $(A)$ with $k\geq 3$, we get the vanishing as follows
$$0=(x_1D_k\circ f)(x_1^2D_2,D_{1k}(x^{\underline{p-1}-2\epsilon_1+2\epsilon_k}))=
-f(x_1^2D_2,D_{1k}(x^{\underline{p-1}-\epsilon_1+\epsilon_k})).$$
On the other hand for the element $D_{12}(x^{\underline{p-1}-\epsilon_1+\epsilon_2})$,
we first use the
following invariance condition
$$0=(x_1D_2\circ f)(x_1^2D_2, D_{12}(x^{\underline{p-1}-2\epsilon_1+2\epsilon_2}))=
[x_1D_2,f(x_1^2D_2,D_{12}(x^{\underline{p-1}-2\epsilon_1+2\epsilon_2}))]+$$
$$-f(x_1^2D_2,D_{12}(x^{\underline{p-1}-\epsilon_1+\epsilon_2})),$$
and then we get the vanishing by means of the following
$$0=(x_1D_2\circ f)(x_1^2D_2,D_{12}(x^{\underline{p-1}-3\epsilon_1+3\epsilon_2}))=
-2f(x_1^2D_2,D_{12}(x^{\underline{p-1}-2\epsilon_1+2\epsilon_2})).$$
Consider now an element $D_{3h}(x^{\underline{a}-2\epsilon_1+\epsilon_3+\epsilon_h})$ of type
$(B)$ and suppose that $a\neq p-2$. Also in this case we get the vanishing using the following
condition
$$0=(x_1D_3\circ f)(x_1^2D_2,D_{3h}(x^{\underline{a}-3\epsilon_1+2\epsilon_3+\epsilon_h}))=
-(a+2)f(x_1^2D_2,D_{3h}(x^{\underline{a}-2\epsilon_1+\epsilon_3+\epsilon_h})).$$
Therefore it remains to consider only the elements of type $(B)$ with $a=p-2$. Define $f(x_1^2D_2, D_{3h}(x^{\underline{p-2}-2\epsilon_1+\epsilon_3+\epsilon_h})):=\gamma_hD_2$
for every $h\neq 3$.
Consider the following invariance conditions for $h\neq 1, 3$:
\begin{equation*}
0=(x_1D_3\circ f)(x_1^2D_2,D_{1h}(x^{\underline{p-2}-2\epsilon_1+\epsilon_3+\epsilon_h}))=
[-\gamma_1 + 4 \gamma_h] D_2 \hspace{0,3cm} \text{ if } p\geq 5, \tag{*}
\end{equation*}
\begin{equation*}
0=(x_1D_3\circ f)(x_1^2D_2,D_{1h}(x^{\underline{1}+\epsilon_1+\epsilon_3+\epsilon_h}))=
\gamma_h D_2 \hspace{2,2cm} \text{ if } p=3, \tag{*'}
\end{equation*}
\begin{equation*}
0=(x_4D_3\circ f)(x_1^2D_2, D_{12}(x^{\underline{p-2}-\epsilon_1+\ep_2+\epsilon_3-\epsilon_4
}))=[-\gamma_1 +3 \gamma_2]D_2 \hspace{0,2cm} \text{ if } n\geq 4. \tag{**}
\end{equation*}

If $n\geq 4$ and $p\geq 5$ then, using $(**)$ and $(*)$ with $h=2$, we get that $\gamma_1=
\gamma_2=0$. Substituting $\gamma_1=0$ in $(*)$, we find $\gamma_h=0$ for every $h$.

If $n\geq 4$ and $p=3$, then from $(*')$, we get the vanishing of $\gamma_h$ for all
$h\neq 1$ and from $(**)$ we get the vanishing of $\gamma_1$.

Finally, if $n=3$ (and $p\geq 5$ by hypothesis) then from $(*)$ we get that $\gamma_1=4\gamma_2$.
We want to prove that if $f\in (E_{\infty}^{1,1})^{S(n)_0}$ then $\gamma_2=0$.
So suppose that $f$ can be lifted to a $S(n)_0$-invariant global cocycle
(which we will continue to call $f$). First of all, by using the $S(n)_0$-invariance
condition $0=(x_2D_3\circ f)(x_1^2D_2,D_{21}(x^{\underline{p-2}-\ep_1+\ep_3}))$, we get
that
$f(x_1^2D_3,D_{21}(x^{\underline{p-2}-\ep_1+\ep_3}))=-5\gamma_2 D_2.$
Using this, we find the following cocycle condition (where we use that $p\geq 5$)
$$0=\d f(x_1^2D_2,x_1^3D_3,D_{12}(x^{\underline{p-2}-3\ep_1+\ep_2+\ep_3}))=
-f(x_1^2D_3,D_{21}(x^{\underline{p-2}-\ep_1+\ep_3}))+$$
$$+f(x_1^2D_2,D_{31}(x^{\underline{p-2}-\ep_1+\ep_3})-5D_{32}(x^{\underline{p-2}-2\ep_1+
\ep_2+\ep_3}))=5\gamma_2 D_2 +4 \gamma_2 D_2 -5\gamma_2 D_2, $$ from
which we deduce that $\gamma_2=0$. \qed

\begin{lemma}\label{S-inv-(0,1)}
Let $d\in \Z_{\geq 0}$.
Then
$C^1\left(S(n)_d,S(n)_{-1}\right)^{S(n)_0}=0.$
\end{lemma}
\dem Obviously a $S(n)_0$-invariant cochain $g \in
C^1(S(n)_d,S(n)_{-1})$ must be homogeneous. Fix $D_i\in S(n)_{-1}$
and let $\phi_i$ be the corresponding weight (hence
$\phi_i=\epsilon_i$ if $i\geq 2$ while
$\phi_1=\sum_{j=2}^n\epsilon_j$).  A base for the space
$S(n)_{\phi_i}$ (which has dimension $(n-1)p$) consists of the
following elements (plus $D_i$):
\begin{align*}
&x_i^{p-1}\otimes T_S, \tag{A}\\
&D_{ij}(x^{\underline{a}+\epsilon_j}) \text{ for } j\neq i \,  \text{ and } 1\leq a\leq p-1
.\tag{B}
\end{align*}
 We have to show that $g$ vanishes on the elements of the above above form.

An element of type $(A)$ must be of the form $x_i^{p-1}D_{jk}(x_jx_k)=x_i^{p-1}(x_jD_j-x_kD_k)$
for some $j,k\neq i$. The vanishing of $g$ on such an element follows from
\begin{equation*}
0=(x_iD_j\circ g)(D_{jk}(x_i^{p-2}x_j^2x_k))=-2 g(x_i^{p-1}D_{jk}(x_jx_k)). \tag{*}
\end{equation*}
Consider now an element $D_{ij}(x^{\underline{a}+\epsilon_j})$ of type $(B)$ and suppose that
$a\neq p-2$. Then we get the vanishing by means of
\begin{equation*}
0=(x_iD_j\circ g)(D_{ji}(x^{\underline{a}+2\epsilon_j-\epsilon_i}))=
(a+2)g(D_{ij}(x^{\underline{a}+\epsilon_j})).\tag{**}
\end{equation*}
Therefore it remains to prove the vanishing for the elements
$D_{ij}(x^{\underline{p-2}+\epsilon_j})$. Put
$g(D_{ij}(x^{\underline{p-2}+\epsilon_j})):=\alpha_j^i D_i$ for $i\neq j$.
Chose three indices $i, j, k$ mutually distinct (which is possible since $n\geq 3$) and consider
the following cocycle condition
$$0=(2x_iD_j\circ g)(D_{jk}(x^{\underline{p-2}+\epsilon_j+\epsilon_k-\epsilon_i}))
=[x_iD_j,g(D_{ik}(x^{\underline{p-2}+\epsilon_k})-D_{ij}(x^{\underline{p-2}+\epsilon_j}))]+$$
\begin{equation*}
+2g(D_{jk}(x^{\underline{p-2}+\epsilon_k}))=(\alpha_j^i-\alpha_k^i+2\alpha_k^j)D_j, \tag{***}
\end{equation*}
where in the first equality we used the relation
$D_{ik}(x^{\underline{p-2}+\epsilon_k})-D_{ij}(x^{\underline{p-2}+\epsilon_j})=
2D_{jk}(x^{\underline{p-2}+\epsilon_j+\epsilon_k-
\epsilon_i}).$
Summing the equation $(***)$ with the one obtained interchanging $k$ with $j$, we get
\begin{equation}
\alpha_k^j+\alpha_j^k=0. \tag{***1}
\end{equation}
Moreover, summing the equation (***) with the analogous one obtained by interchanging $i$ with $j$
and using the antisymmetric property (***1), we obtain
\begin{equation*}
\alpha_k^i+\alpha_k^j=0. \tag{***2}
\end{equation*}
Finally, using equations $(***1)$ and $(***2)$, we get
$\alpha_j^i=-\alpha_j^k=\alpha_k^j$ and $\alpha_k^i=-\alpha_k^j$.
Substituting into the equation (***), we find $4\alpha_k^j=0$. \qed

\begin{lemma}\label{S-H^1-(0,1)}
Let $d\in \Z_{\geq 0}$. Then
$$H^1\left(S(n)_0,C^1\left(S(n)_d,S(n)_{-1}\right)\right)=
\begin{cases}
\oplus_{i=1}^n F \cdot \langle  \overline{{\rm Sq}(D_i)} \rangle &\text{ if } d=p-1,\\
 0 &\text{ otherwise, }
\end{cases}$$
where $\overline{{\rm Sq}(D_i)}$ denotes the restriction of ${\rm Sq}(D_i)$ to $S(n)_0\times
S(n)_{p-1}$.
\end{lemma}
\dem First of all, observe that the computations made at the
beginning of subsection 4.2 show that the above cocycles
$\overline{{\rm Sq}(D_i)}$ are independent modulo coboundaries.
Consider a cocycle $f\in \bigoplus_{d\geq 0}
Z^1(S(n)_0,C^1(S(n)_d,S(n)_{-1}))$. Since the maximal torus $T_S$ is
contained in $S(n)_0$, we can assume that $f$ is homogeneous.
Exactly as in the proof of lemma \ref{W-H^1-(0,1)}, one can show,
using the above lemma \ref{S-inv-(0,1)}, that the restriction of $f$
to the maximal torus $T_S$ is zero. Therefore, by homogeneity, the
cocycle $f$ can take only the following non-zero values (with $1\leq
i,j,k\leq n$ mutually distinct):
\begin{align*}
&f_{x_iD_j}(E)\subset\langle D_j\rangle  \text{ if } E=
\begin{sis}
&x_i^{p-1}D_{kh}(x_kx_h)  \text{ for } h\neq i,k,  &\hspace{1cm}  \text{ (1A) } \\
&D_{ih}(x^{\underline{a}+\epsilon_h}) \text{ for } a\neq 0, p-2 \text{ and } h\neq i,
 &\text{ (1B) } \\
& D_{jh}(x^{\underline{p-2}-\epsilon_i+\epsilon_j+\epsilon_h}) \text{ for } h\neq j,
 & \text{ (1C) }
\end{sis}\\
&f_{x_iD_j}(E)\subset\langle D_i\rangle \text{ if } E=
\begin{sis}
&D_{jh}(x_i^{p-2}x_j^2x_h)  \text{ for } h\neq  j, &  \text{ (2A) }\\
&D_{jh}(x^{\underline{a}+2\epsilon_j-2\epsilon_i+\epsilon_h}) \text{
for } a\neq 0, p-2, \:
h\neq j,& \text{ (2B) }\\
& D_{kh}(x^{\underline{p-2}-2\epsilon_i+\epsilon_j+\epsilon_k+\epsilon_h}) \text{ for } h\neq k,
& \text{ (2C) }
\end{sis}\\
&f_{x_iD_j}(E)\subset\langle D_k\rangle \text{ if } E=
\begin{sis}
& D_{jh}(x^{-\epsilon_i+2\epsilon_j-\epsilon_k+\epsilon_h})  \text{ for } h\neq j, & \text{ (3A) }\\
&
D_{jh}(x^{\underline{a}-\epsilon_i+2\epsilon_j-\epsilon_k+\epsilon_h})
\text{ for } a\neq 0, p-2, \:
 h\neq j, &\text{ (3B) }\\
& D_{kh}(x^{\underline{p-2}-\epsilon_i+\epsilon_j+\epsilon_h}) \text{ for } h\neq k. &\text{ (3C) }
\end{sis}
\end{align*}
We want to show that we can modify $f$, by adding coboundaries and
the cocycles $\overline{\Sq(D_i)}$, in such a way that it vanishes
on the above elements. We divide the proof in several steps
according to the elements of the above list.

$\fbox{(2A)}$  For every index $i$, we choose an index $j\neq i$ and
we modify $f$, by adding a multiple of $\overline{{\rm Sq}(D_i)}$,
in such a way that $f_{x_iD_j}(D_{ji}(x_i^{p-1}x_j^2))=0$ (see
equation (\ref{S-ind-coc})). Moreover, by adding a coboundary $\d
g$, we can further modify $f$ in such a way that
$f_{x_iD_j}(D_{jk}(x_i^{p-2}x_j^2x_k))=0$ for every $k\neq i, j$
(see equation (*) of lemma \ref{S-inv-(0,1)}). Therefore we get
the required vanishing for the chosen index $j$. Using this, we
obtain the following cocycle condition (for every $k\neq i,j$ and
$h\neq j$):
$$0=\d f_{(x_iD_j,x_iD_k)}(D_{jh}(x^{-3\epsilon_i+\epsilon_j+\epsilon_k+\epsilon_h}))=
-2f_{x_iD_k}(D_{jh}(x^{-2\epsilon_i+\epsilon_j+\epsilon_k+\epsilon_h})),
$$
from which we get the required vanishing, using (for $h\neq k$) the
transformation rule
$D_{kh}(x_i^{p-2}x_k^2x_h)=2D_{jh}(x^{-2\epsilon_i+\epsilon_j+\epsilon_k+\epsilon_h})-
D_{jk}(x^{-2\epsilon_i+\epsilon_j+2\epsilon_k})$.

$\fbox{(3A)}$
 If $p\geq 5$ then we get the required vanishing by means of the
following condition, where we used the vanishing of the elements of
type $(2A)$:
$$0=\d f_{(x_kD_j,x_iD_j)}(D_{jh}(x^{-2\epsilon_k+3\epsilon_j-\epsilon_k+\epsilon_h}))=
-3
f_{x_iD_j}(D_{jh}(x^{-\epsilon_i+2\epsilon_j-\epsilon_k+\epsilon_h})).$$
If $p=3$ a little extra-work is necessary and we have to consider
the following three conditions according to the three cases $h\neq
i,k$, $h=i$ and $h=k$ respectively:
$$0=\d f_{(x_kD_h,x_iD_j)}(x_i^2x_k^2D_{jh}(x_j^2x_h))=[x_kD_h,f_{x_iD_j}(x_i^2x_k^2D_{jh}(x_j^2x_h))],$$
$$0=\d f_{x_iD_j,x_kD_j)}(x_kx_j^2D_i)=-f_{x_iD_j}(2x_k^2x_jD_i)=f_{x_iD_j}(x_k^2D_{ji}(x_j^2)),$$
$$0=\d f_{(x_iD_j,x_kD_i)}(x_i^2D_{jk}(x_j^2))=-[x_kD_i,f_{x_iD_j}(x_i^2D_{jk}(x_j^2))],$$
where in the last condition we used the  first two vanishing.

$\fbox{(1A)}$
Using the vanishing of (2A) and (3A), we get
$$0=\d f_{(x_iD_j,x_iD_k)}(x_i^{p-2}D_{kh}(x_k^2x_h))=
2f_{x_iD_j}(x_i^{p-1}D_{kh}(x_kx_h)).$$

$\fbox{(2B)}$
 Fix an integer $1\leq a \leq p-1$ different from $p-2$
and define\\
$f_{x_iD_j}(D_{jh}(x^{\underline{a}+2\epsilon_j-2\epsilon_i+\epsilon_h}))=
\gamma^i_{jh} D_i$ for every $j\neq i, h$.  By adding a coboundary
$\d g$, we can modify $f$ in such a way that $\gamma^i_{ji}=0$ for
every $j\neq i$ (see equation (**) of lemma \ref{S-inv-(0,1)}).
Consider first the cocycle condition (for $i,j,k$ mutually distinct)
$$0=\d f_{(x_iD_j,x_iD_k)}(D_{ij}(x^{\underline{a}-2\epsilon_i+2\epsilon_j+\epsilon_k})=
(a+2)f_{x_iD_j}(D_{ji}(x^{\underline{a}-\epsilon_i+\epsilon_j+\epsilon_k}))+$$
$$-(a+2)f_{x_iD_j}(D_{ki}(x^{\underline{a}-\epsilon_i+\epsilon_j+\epsilon_k}))-
(a-2)\gamma^i_{jk}D_i.$$ By considering the analogous condition
obtained by interchanging $j$ with $k$ together with the
transformation rule
$$(a+2)D_{ji}(x^{\underline{a}-\epsilon_i+\epsilon_j
+\epsilon_k})=(a+1)D_{ki}(x^{\underline{a}-\epsilon_i+2\epsilon_k})-(a-1)D_{kj}(x^{\underline{a}
-2\epsilon_i+2\epsilon_k+\epsilon_j}),$$ we get the relation
$(1-a)\gamma^i_{kj}+\gamma^i_{jk}=0.$
 Next consider the other
cocycle condition
$$0=\d f_{(x_iD_j,x_iD_k)}(D_{jk}(x^{\underline{a}-3\epsilon_i+2\epsilon_j+2\epsilon_k}))=
(a+2)(\gamma^i_{kj}+\gamma^i_{jk})D_i.$$
Since ${\rm det}
\begin{pmatrix}
1-a& 1\\
a+2& a+2\\
\end{pmatrix}
=-a(a+2)\neq 0$, putting together these two relations we get that
$\gamma^i_{jk}=0$.

$\fbox{(3B)}$
 If $a\neq p-3$ then, using the vanishing of the elements of type (2B), we get
$$0=\d f_{(x_kD_j,x_ID_j)}(D_{jh}(x^{\underline{a}-2\ep_k+3\ep_j-\ep_i+\ep_h}))=
-(a+3)f_{x_iD_j}(D_{jh}(x^{\underline{a}-\ep_i+2\ep_j-\ep_k+\ep_h})).$$
If $a=p-3$ (and hence $p\geq 5$) we use the following condition (again by the
vanishing of (2B))
$$0=\d f_{(x_kD_i,x_iD_j)}(D_{ih}(x^{\underline{p-3}-2\ep_k+\ep_i+\ep_j+\ep_h}))=
2f_{x_iD_j}(D_{ih}(x^{\underline{p-3}-\ep_k+\ep_j+\ep_h})),$$
together with the transformation rule (if $h\neq i,j$)
$$3D_{jh}(x^{\underline{p-3}-\ep_i+2\ep_j-\ep_k+\ep_h})=-(2+\delta_{hk})
D_{ij}(x^{\underline{p-3}-\ep_k+2\ep_j})+D_{ih}(x^{\underline{p-3}-\ep_k+\ep_j+\ep_h}).$$

$\fbox{(1B)}$
 Take indices $r\neq i,j$ and $s\neq r$ and consider
the following condition (using the vanishing of (2B) and (3B))
$$0=\d f_{(x_iD_j,x_iD_r)}(D_{rs}(x^{\underline{a}+\ep_r-\ep_i+\ep_s}))=(a+2)
f_{x_iD_j}(D_{rs}(x^{\underline{a}-\ep_i+\ep_r+\ep_s})).$$ By taking
$r=h$ and $s=i$, we get the required vanishing if $h\neq j$. If
$h=j$ and $a\neq p-1$, we use the transformation rule
$$(a+1)D_{ij}(x^{\underline{a}+\ep_j})=(a+1)^2D_{ri}(x^{\underline{a}+\ep_r})-a(a+1)
D_{rj}(x^{\underline{a}-\ep_i+\ep_r+\ep_j}).$$ If $h=j$ and $a=p-1$
we use the following condition (by the vanishing of (2B) and (3B))
$$0=\d f_{(x_kD_j,x_iD_j)}(D_{ji}(x^{\underline{p-1}-\ep_k+2\ep_j}))=f_{x_iD_j}
(D_{ij}(x^{\underline{p-1}+\ep_j})).$$

$\fbox{(1C)}$
 We define
$f_{x_iD_j}(D_{jk}(x^{\underline{p-2}-\epsilon_i+\epsilon_j+\epsilon_k}))=
\beta^i_{jk}D_j$ for every $j\neq i,k$ (but possibly $i=k$). The
space of all such cochains has dimension $n(n-1)^2$. Using the
notations of lemma \ref{S-inv-(0,1)}, the subspace of coboundaries
is formed by the $\beta_{jk}^i$ such that there exist
$\{\alpha_j^i\: :\: i\neq j\}$ with the property that
$2\beta_{jk}^i=\alpha_j^i- \alpha_k^i+2\alpha_k^j$ (see equation
(***) of lemma \ref{S-inv-(0,1)}). Moreover in the above quoted
lemma, we prove that different values of $\alpha_j^i$ give rise to
different values of $\beta_{jk}^i$. Hence the dimension of the
subspace of coboundaries is $n(n-1)$. Therefore, in order to prove
the vanishing of the elements of type $(1C)$, it will be enough to
exhibit $n(n-1)(n-2)$ linearly independent relations among the
coefficients $\beta_{jk}^i$.

Fix three integers $i, j, k$ mutually distinct and consider
the following cocycle condition
$$0=\d f_{(x_iD_j,x_jD_k)}(D_{jk}(x^{\underline{p-2}+\epsilon_j+\epsilon_k-\epsilon_i}))=
(-\beta_{kj}^j+\beta_{jk}^i+\beta_{kj}^i)D_k.$$ We get first of all
that the $\beta$'s with two coincident indices are determined by
those with three different indices and this give $n(n-1)$ linearly
independent relations. Moreover we deduce also that for any $k\neq
j$ the value of the sum $\beta_{jk}^i+\beta_{kj}^i$ is independent
of $i$ and this give $n(n-1)(n-3)$ linearly independent relations.
Since the two types of relations are also independent one of the
other, the total number of independent relations we get is
$n(n-1)(n-2)$, as required.

 $\fbox{(3C)}$
 Using the vanishing of (1C), we get
$$0=\d f_{(x_iD_j,x_jD_k)}(D_{kh}(x^{\underline{p-2}+\ep_k+\ep_h-\ep_i}))=
-f_{x_iD_j}(D_{kh}(x^{\underline{p-2}+\ep_h-\ep_i+\ep_j})).$$

$\fbox{(2C)}$
 Using the vanishing of (1C) and (3C), we compute
\begin{equation*}
\begin{split}
& 0=\d f_{(x_iD_k,x_iD_j)}(D_{kh}(x^{\underline{p-2}-2\ep_i+\ep_j+
\ep_k+\ep_h}))= \\
&=[x_iD_k,f_{x_iD_j}(D_{kh}(x^{\underline{p-2}-2\ep_i+\ep_j+\ep_k+
\ep_h}))].
\end{split} 
\end{equation*}
\qed


\begin{thebibliography}{TESI}


\bibitem[BW88]{BW} R.E. Block, R.L. Wilson, Classification of the restricted simple
Lie algebras,  J. Algebra  114  (1988) 115--259.


\bibitem[CHE05]{Che1} N. G. Chebochko, Deformations of classical Lie algebras with a
homogeneous root system in characteristic two I (Russian), Mat. Sb.
196  (2005) 125--156. English translation:  Sb. Math.  196 (2005)
1371--1402.

\bibitem[CK00]{Che2} N. G. Chebochko, M. I. Kuznetsov, Deformations of classical Lie
algebras (Russian), Mat. Sb. 191 (2000) 69--88. English translation:
Sb. Math. 191 (2000) 1171--1190.

\bibitem[CKK00]{Che3} N. G. Chebochko, S. A. Kirillov, M. I.
Kuznetsov, Deformations of a Lie algebra of type $G\sb 2$ of
characteristic three (Russian), Izv. Vyssh. Uchebn. Zaved. Mat. no.
3 (2000) 33--38. English translation: Russian Math. (Iz. VUZ) 44
(2000) 31--36.


\bibitem[CEL70]{CEL} M. Ju. Celousov, Derivations of Lie algebras of
Cartan type (Russian), Izv. Vys\v s. U\v cebn. Zaved. Matematika 98
(1970) 126--134.


\bibitem[CE48]{CE} C. Chevalley, S. Eilenberg, Cohomology theory of Lie groups
and Lie algebras, Trans. Amer. Math. Soc. 63 (1948) 85--124.

\bibitem[DG70]{DG} M. Demazure, P. Gabriel, Groupes alg\'ebriques,
Tome I: G\'eom\'etrie alg\'ebrique, g\'en\'eralit\'es, groupes
commutatifs (French), Masson and Cie Editeur, North-Holland
Publishing Co., Amsterdam, 1970.

\bibitem[DK78]{DK} A. S. D\v zumadildaev, A. I. Kostrikin,
Deformations of the Lie algebra $W\sb{1}(m)$ (Russian), Algebra,
number theory and their applications,  Trudy Mat. Inst. Steklov. 148
(1978) 141--155.

\bibitem[DZU80]{DZ1} A. S. D\v zumadildaev, Deformations of general Lie algebras
of Cartan type (Russian),  Dokl. Akad. Nauk SSSR  251  (1980)
1289--1292. English translation: Soviet Math. Dokl. 21 (1980)
605--609.

\bibitem[DZU81]{DZ2} A. S. D\v zumadildaev, Relative cohomology and deformations
of the Lie algebras of Cartan types (Russian),  Dokl. Akad. Nauk
SSSR 257  (1981) 1044--1048. English translation: Soviet Math. Dokl.
23 (1981) 398--402.

\bibitem[DZU89]{DZ3} A. S. D\v zumadildaev, Deformations of the Lie algebras
$W\sb n(m)$ (Russian),  Mat. Sb.  180  (1989) 168--186. English
translation: Math. USSR-Sb.  66  (1990) 169--187.

\bibitem[FS88]{FS} R. Farnsteiner, H. Strade, Modular Lie algebras
and their representation, Monographs and textbooks in pure and
applied mathematics vol. 116., Dekker, New York, 1988.

\bibitem[GER64]{GER1} M. Gerstenhaber, On the deformation of rings and
algebras, Ann. of Math. 79  (1964) 59--103.

\bibitem[HS97]{HiSt} P. J. Hilton, U. A. Stammbach, A course in homological
algebra, Graduate Texts in Mathematics 4, Springer-Verlag, New York,
1997.


\bibitem[HS53]{HS} G. Hochschild, J.-P. Serre, Cohomology of Lie
algebras, Ann. of Math. 57 (1953) 591--603.

\bibitem[KS66]{KS} A. I. Kostrikin, I. R. Shafarevich, Cartan's pseudogroups
and the $p$-algebras of Lie (Russian), Dokl. Akad. Nauk SSSR  168
(1966) 740--742. English translation: Soviet Math. Dokl. 7 (1966)
715--718.

\bibitem[MEL80]{MEL} G. M. Melikian, Simple Lie algebras of characteristic $5$
(Russian),  Uspekhi Mat. Nauk  35  (1980) 203--204.

\bibitem[PS97]{PS1} A. Premet, H. Strade, Simple Lie algebras of small characteristic
I: Sandwich elements,  J. Algebra  189  (1997) 419--480.

\bibitem[PS99]{PS2} A. Premet, H. Strade, Simple Lie algebras of small characteristic
II: Exceptional roots,  J. Algebra  216  (1999) 190--301.

\bibitem[PS01]{PS3} A. Premet, H. Strade, Simple Lie algebras of small characteristic
III: The toral rank 2 case,   J. Algebra  242  (2001) 236--337.

\bibitem[RUD71]{RUD} A. N. Rudakov, Deformations of simple Lie algebras
(Russian), Izv. Akad. Nauk SSSR Ser. Mat. 35 (1971) 1113--1119.


\bibitem[SEL67]{SEL} G. B. Seligman, Modular Lie algebras, Ergebnisse der Mathematik
und ihrer Grenzgebiete Band 40, Springer-Verlag, New York, 1967.


\bibitem[STR89]{STR1} H. Strade, The Classification of the Simple Modular Lie Algebras
I: Determination of the two-sections, Ann. of Math. 130 (1989)
643--677.

\bibitem[STR92]{STR2} H. Strade, The Classification of the Simple Modular Lie Algebras
II: The Toral Structure, J. Algebra 151 (1992) 425--475.

\bibitem[STR91]{STR3} H. Strade, The Classification of the Simple Modular Lie Algebras
III: Solution of the Classical Case, Ann. of Math. 133 (1991)
 577--604.

\bibitem[STR93]{STR4} H. Strade, The Classification of the Simple Modular Lie Algebras
IV: Determining the Associated Graded Algebra, Ann. of Math. 138
(1993) 1--59.

\bibitem[STR94]{STR5} H. Strade, The Classification of the Simple Modular Lie Algebras
V: Algebras with Hamiltonian Two-sections, Abh. Math. Sem. Univ.
Hamburg 64 (1994) 167--202.

\bibitem[STR98]{STR6} H. Strade, The classification of the simple modular Lie algebras
VI: Solving the final case,  Trans. Amer. Math. Soc.  350  (1998)
2553--2628.

\bibitem[STR04]{STR} H. Strade, Simple Lie algebras over fields of positive characteristic I:
Structure theory, De Gruyter Expositions in Mathematics 38, Walter
de Gruyter, Berlin, 2004.

\bibitem[VIV2]{VIV2} F. Viviani, Infinitesimal deformations of restricted simple 
Lie algebras II, submitted (available at arXiv:math/0702499).

\bibitem[VIV3]{VIV3} F. Viviani, Deformations of the restricted Melikian algebra,
to appear on Comm. Algebra.

\bibitem[VIV4]{VIV4} F. Viviani, Deformations of simple finite group schemes,
submitted (available at arXiv:0705.0821).

\bibitem[SW91]{SW} R. L. Wilson, H. Strade, Classification of Simple Lie Algebras over
Algebraically Closed Fields of Prime Characteristic, Bull. Amer.
Math. Soc. 24 (1991) 357--362.





\end{thebibliography}
\end{document}